\documentclass[journal,twoside]{IEEEtran}

\usepackage{cite}
\usepackage[cmex10]{amsmath,mathtools}
\usepackage{amsfonts}
\usepackage{amssymb}
\usepackage{textcomp}
\usepackage{mathrsfs}
\usepackage{inconsolata}
\usepackage[caption=false,font=footnotesize]{subfig}
\usepackage{fixltx2e}
\usepackage{url}
\usepackage[implicit=false]{hyperref}
\usepackage{siunitx}
\usepackage{xspace}
\usepackage{accents}
\usepackage{ifthen}
\usepackage{blindtext}
\usepackage{bm}
\usepackage[multiple]{footmisc}
\MakeRobust{\underaccent}

\newtheorem{theorem}{Theorem}
\newtheorem{lemma}{Lemma}
\newtheorem{proposition}{Proposition}

\newtheorem{definition}{Definition}
\newtheorem{requirement}[definition]{Definition}

\DeclareMathOperator{\real}{Re}
\DeclareMathOperator{\imag}{Im}
\DeclareMathOperator{\trace}{tr}
\DeclareMathOperator{\rank}{rank}

\DeclareMathOperator*{\minimize}{minimize}

\DeclareMathOperator*{\subjectto}{subject~to}

\DeclareMathOperator{\psd}{\hskip1pt\succeq\hskip1pt}

\newcommand{\abs}[1]{{\lvert#1\rvert}}
\newcommand{\norm}[1]{{\lVert#1\rVert}}
\newcommand{\vecmtx}[1]{{\bm{#1}}}
\newcommand{\mthcal}[1]{{\mathcal{#1}}}

\newcommand{\conj}{^\ast}
\DeclareMathOperator{\transpose}{T}
\DeclareMathOperator{\hermitian}{H}
\newcommand{\tran}{^{\transpose}}
\newcommand{\herm}{^{\hermitian}}

\newcommand{\nR}{\mathbb{R}}
\newcommand{\nRnn}{\mathbb{R}_{+}}
\newcommand{\nRp}{\mathbb{R}_{++}}
\newcommand{\nZ}{\mathbb{Z}}
\newcommand{\nN}{\mathbb{N}}
\newcommand{\nC}{\mathbb{C}}
\newcommand{\nS}{\mathbb{S}}
\newcommand{\nSp}{\tilde{\mathbb{S}}}
\newcommand{\aSrc}[1]{\expandafter\hat#1}
\newcommand{\aDst}[1]{\expandafter\check#1}
\newcommand{\aFwd}[1]{\expandafter\hat#1}
\newcommand{\aBwd}[1]{\expandafter\check#1}
\newcommand{\aLB}[1]{\expandafter\underaccent{\bar}#1}
\newcommand{\aUB}[1]{\expandafter\bar#1}

\newcommand{\aSrcUB}[1]{\bar{\hat{#1}}}

\newcommand{\aDstUB}[1]{\bar{\check{#1}}}
\newcommand{\aAC}[1]{\expandafter\tilde#1}
\newcommand{\aDC}[1]{\expandafter\bar#1}
\newcommand{\aSh}[1]{\expandafter\tilde#1}
\newcommand{\aBar}[1]{\expandafter\bar#1}
\newcommand{\aUBar}[1]{\expandafter\underaccent{\bar}#1}

\newcommand{\aIj}[1]{#1^{\mathrm{in}}}
\newcommand{\aP}[1]{#1^{\mathrm{p}}}
\newcommand{\aQ}[1]{#1^{\mathrm{q}}}
\newcommand{\aOpt}[1]{#1^{\star}}
\newcommand{\aUndirected}[1]{\expandafter\underaccent{\bar}#1}
\newcommand{\aEQ}[1]{\expandafter\bar#1}
\newcommand{\aIE}[1]{\expandafter\check#1}

\newcommand{\sB}{{\mthcal{B}}}
\newcommand{\sC}{{\mthcal{C}}}

\newcommand{\sE}{{\mthcal{E}}}
\newcommand{\sF}{{\mthcal{F}}}

\newcommand{\sI}{{\mthcal{I}}}
\newcommand{\sJ}{{\mthcal{J}}}

\newcommand{\sM}{{\mthcal{M}}}
\newcommand{\sN}{{\mthcal{N}}}

\newcommand{\sS}{{\mthcal{S}}}

\newcommand{\sU}{{\mthcal{U}}}
\newcommand{\sV}{{\mthcal{V}}}

\newcommand{\sX}{{\mthcal{X}}}

\newcommand{\sVAC}{{\aAC{\sV}}}
\newcommand{\sVDC}{{\aDC{\sV}}}
\newcommand{\sEAC}{{\aAC{\sE}}}
\newcommand{\sEDC}{{\aDC{\sE}}}
\newcommand{\sBE}{\sB_\sE}
\newcommand{\sBC}{\sB_\sC}
\newcommand{\sBI}{\sB_\sI}

\newcommand{\cC}{{\abs{\sC}}}

\newcommand{\cE}{{\abs{\sE}}}

\newcommand{\cI}{{\abs{\sI}}}
\newcommand{\cJ}{{\abs{\sJ}}}

\newcommand{\cM}{{\abs{\sM}}}

\newcommand{\cV}{{\abs{\sV}}}

\newcommand{\cEAC}{{\abs{\sEAC}}}

\newcommand{\gG}{{\mthcal{G}}}
\newcommand{\eSrc}{{\aSrc{\epsilon}}}
\newcommand{\eDst}{{\aDst{\epsilon}}}
\newcommand{\cSrc}{{\aSrc{\gamma}}}
\newcommand{\cDst}{{\aDst{\gamma}}}
\newcommand{\nSrc}{{\aSrc{n}}}

\newcommand{\fC}{C}

\newcommand{\fLoss}{L}

\newcommand{\fProj}{\mathscr{P}}

\newcommand{\mA}{{\vecmtx{A}}}
\newcommand{\mB}{{\vecmtx{B}}}
\newcommand{\mC}{{\vecmtx{C}}}

\newcommand{\mH}{{\vecmtx{H}}}

\newcommand{\mL}{{\vecmtx{L}}}
\newcommand{\mM}{{\vecmtx{M}}}

\newcommand{\mP}{{\vecmtx{P}}}
\newcommand{\mQ}{{\vecmtx{Q}}}

\newcommand{\mS}{{\vecmtx{S}}}

\newcommand{\mV}{{\vecmtx{V}}}

\newcommand{\mY}{{\vecmtx{Y}}}

\newcommand{\mZero}{{\vecmtx{0}}}

\newcommand{\mISrc}{{\hat{\skew{-3}\bm{I}}\hskip-2.25pt}}
\newcommand{\mIDst}{{\check{\skew{-3}\bm{I}}\hskip-2.25pt}}
\newcommand{\mAUB}{{\bar{\skew{-6}\bm{A}}\hskip-3.5pt}}
\newcommand{\mVp}{\mV}

\newcommand{\vc}{{\vecmtx{c}}}

\newcommand{\ve}{{\vecmtx{e}}}
\newcommand{\vf}{{\vecmtx{f}}}

\newcommand{\vh}{{\vecmtx{h}}}
\newcommand{\vi}{{\vecmtx{i}}}

\newcommand{\vl}{{\vecmtx{l}}}

\newcommand{\vp}{{\vecmtx{p}}}
\newcommand{\vq}{{\vecmtx{q}}}

\newcommand{\vs}{{\vecmtx{s}}}

\newcommand{\vv}{{\vecmtx{v}}}
\newcommand{\vw}{{\vecmtx{w}}}
\newcommand{\vx}{{\vecmtx{x}}}

\newcommand{\viSrc}{{\aSrc{\skew{-2.5}\vecmtx{i}}}}
\newcommand{\viDst}{{\aDst{\skew{-2.5}\vecmtx{i}}}}
\newcommand{\hynet}{\emph{hynet}\xspace}
\newcommand{\matpower}{\textsc{Matpower}\xspace}
\newcommand{\pypower}{\textsc{Pypower}\xspace}
\newcommand{\powermodels}{PowerModels\xspace}
\newcommand{\psat}{PSAT\xspace}
\newcommand{\pandapower}{\texttt{pandapower}\xspace}

\newcommand{\matacdc}{MATACDC\xspace}
\newcommand{\iu}{{\mkern1.5mu\text{\bfseries i}\mkern1.5mu}}

\begin{document}

\title{\emph{hynet:} An Optimal Power Flow Framework for Hybrid AC/DC Power Systems}

\author{Matthias~Hotz,~\IEEEmembership{Student Member,~IEEE,}~and~Wolfgang~Utschick,~\IEEEmembership{Senior Member,~IEEE}%
\thanks{The authors are with the Department of Electrical and Computer Engineering, Technische Universit\"at M\"unchen, Munich D-80333, Germany (e-mail: matthias.hotz@tum.de, utschick@tum.de).}%
}

\markboth{Hotz and Utschick: hynet -- An Optimal Power Flow Framework for Hybrid AC/DC Power Systems}{Hotz and Utschick: hynet -- An Optimal Power Flow Framework for Hybrid AC/DC Power Systems}

\maketitle

\begin{abstract}
High-voltage direct current (HVDC) systems are increasingly incorporated into today's AC power grids, necessitating optimal power flow (OPF) tools for the analysis, planning, and operation of such hybrid systems. To this end, we introduce \hynet, a Python-based open-source OPF framework for hybrid AC/DC grids with point-to-point and radial multi-terminal HVDC systems. \hynet's design promotes ease of use and extensibility, which is supported by the particular mathematical model and software design presented in this paper. The system model features a unified representation of AC and DC subgrids as well as a concise and flexible converter model, which enable the compact description of a hybrid AC/DC power system and its OPF problem. To support convex relaxation based OPF solution techniques, a state space relaxation is introduced to obtain a unified OPF formulation that is analogous to the OPF of AC power systems. This enables the direct generalization of relaxation-related results for AC grids to hybrid AC/DC grids, which is shown for the semidefinite and second-order cone relaxation as well as associated results on exactness and locational marginal prices. Finally, \hynet's object-oriented software design is discussed, which provides extensibility via inheritance and standard design patterns, and its robust and competitive performance is illustrated with case studies.
\end{abstract}


\section{Introduction}
\label{sec:introduction}

\IEEEPARstart{T}{o} counteract the climate change, many countries consider a decarbonization of the energy sector, especially via a transition of electricity generation based on fossil fuels toward renewable energy sources (RES)~\cite{European-Commission2017a,Arcia-Garibaldi2018a}. This transition introduces an increasingly distributed and fluctuating energy production, which generally necessitates additional transmission capacity as well as stronger interconnections of regional and national grids to balance and smooth the variability of RES-based generation~\cite{Arcia-Garibaldi2018a,IEA2016a}. In this regard, high-voltage direct current (HVDC) systems are considered as a key technology due to their advantages in long-distance, underground, and submarine transmission as well as their ability to connect asynchronous grids and, in case of voltage source converter (VSC) HVDC systems, to provide flexible power flow control and reactive power compensation~\cite{Arcia-Garibaldi2018a,IEA2016a}. Already today, a large number of point-to-point HVDC (P2P-HVDC) systems and several multi-terminal HVDC (MT-HVDC) systems are installed and many are planned~\cite{Arcia-Garibaldi2018a,Buigues2017a}. With the above developments, this trend is destined to continue, leading to large-scale hybrid AC/DC power systems.

Due to the importance of the optimal power flow (OPF) in operational and grid expansion planning as well as techno-economic studies, these structural changes necessitate an OPF framework for large-scale hybrid AC/DC power systems. OPF denotes the optimization problem of identifying the cost-optimal allocation of generation resources and the corresponding system state to serve a given load, while satisfying all boundary conditions of the grid. It involves a large number of optimization variables, system constraints, and, to accurately capture the physics, the power flow equations based on Kirchhoff's laws, which render the problem inherently nonconvex and challenging to solve. As the OPF problem specification and solution is rather involved, a software framework for this task is desired. Furthermore, for transparency, reproducibility, and flexible adoption in research, it should be available as open-source software. Several open-source software packages for OPF computation have already been published, including the established toolboxes \matpower~\cite{Zimmerman2011a,Zimmerman2016a} (and its Python-port \pypower~\cite{Lincoln2017a}) and \psat~\cite{Milano2005a,Milano2016a} as well as the recently released \powermodels~\cite{Coffrin2018a,Coffrin2018b} and \pandapower~\cite{Thurner2018a,Thurner2018b}. While \psat is targeted at small to medium-sized systems, \matpower, \powermodels, and \pandapower also support large-scale systems, but they are limited to a simple model of P2P-HVDC systems and do not support MT-HVDC systems. On the contrary, the open-source software \matacdc~\cite{Beerten2015a} features an elaborate model for hybrid AC/DC grids, but it is limited to (sequential) power flow computations.

In the literature, several OPF formulations for hybrid AC/DC grids with MT-HVDC systems were proposed, e.g.~\cite{Wiget2012a,Cao2013a,Feng2014a,Zhao2017a,Baradar2013a,Bahrami2017a}, which extend the standard AC model with DC buses, DC lines and cables, and VSCs. While the modeling of DC grids, where DC lines and cables are considered via their equivalent series resistance, resembles that of AC grids and is rather straightforward, the modeling of VSCs is more intricate. Elaborate models for VSC stations comprise a transformer, filter, phase reactor, and converter. While the former three can be considered via equivalent $\pi$-models, the converter is modeled as a lossy transfer of active power between the AC and DC side of the converter with the provision of reactive power on the AC side, cf.~\cite{Zhao2017a,Beerten2012a,Cigre2014b}. Many recent works model the converter losses as a quadratic function of the converter current~\cite{Wiget2012a,Cao2013a,Feng2014a,Bahrami2017a,Beerten2012a,Baradar2013b} which traces back to the analysis in~\cite{Daelemans2008a}, with the most advanced formulation in~\cite{Zhao2017a} using an individual parameterization for the rectifier and inverter mode. In terms of operating limits, the most elaborate characterizations include a converter current limit as well as constraints that restrict the provision of reactive power based on the coupling of the AC- and DC-side voltage, cf.~\cite{Zhao2017a,Feng2014a,Beerten2012a}. As a consequence of the significant increase in model complexity compared to AC systems, the mathematical formulation and parameterization of the OPF problem for hybrid AC/DC power systems is intricate and extensive.

The OPF problem, even for AC-only grids, is notoriously hard to solve and a recent technique to improve the computational tractability of its globally optimal solution is \emph{convex relaxation}. Instead of the traditional approach of simplifying the system model, e.g., to the widely used ``DC power flow''~\cite{Wood1996a}, convex relaxation simplifies the structure of the feasible set. In the context of (optimal) power flow, the first second-order cone and semidefinite relaxation is attributed to Jabr~\cite{Jabr2006a} and Bai \emph{et al.}~\cite{Bai2008a}, respectively, and was followed by a vast number of works, see e.g.~\cite{Low2014a,Taylor2015a,Molzahn2015a,Coffrin2016a,Kocuk2016b,Bingane2018a} and the references therein. However, besides the huge body of literature on convex relaxation for AC grids, only some works consider hybrid AC/DC grids with P2P-~\cite{Hotz2016a,Hotz2018a} and MT-HVDC systems~\cite{Baradar2013a,Bahrami2017a,Bahrami2018a,Venzke2019a}.

\subsection{Contributions}

To support the study of hybrid AC/DC power systems, we developed \hynet~\cite{Hotz2018b}, an open-source OPF framework for hybrid AC/DC grids with P2P- and radial MT-HVDC systems. In its design, we emphasized ease of use and extensibility to facilitate its effortless adoption in research and education. To this end, we developed a solid and rigorous yet compact and easily parameterizable system model as well as a specialized OPF formulation that simplifies the study of convex relaxations. This mathematical foundation was embedded in a clearly structured object-oriented software design to create \hynet, which was written in the popular high-level open-source programming language Python~\cite{Python-Core-Team2018a} that is freely available for all major platforms. While the documentation of the software framework is provided on the \hynet project website~\cite{Hotz2018b}, this paper focuses on \hynet's mathematical foundation as well as its relation to the software and comprises the following main contributions:
\begin{itemize}
	\item \emph{System model for hybrid AC/DC power systems:}
		By introducing a new network topology model as well as a set of definitions, we generalize the electrical model for AC grids in~\cite{Hotz2016a} and~\cite{Hotz2018a} to AC and DC subgrids and develop a \emph{unified} AC/DC subgrid model. This is complemented by a new converter model that balances modeling depth and complexity to offer adequate accuracy based on a concise and flexible parameterization. Jointly, they enable a compact description of hybrid AC/DC systems. The presented rigorous model formulation does not only provide a solid mathematical foundation for \hynet, but it also includes a manifold of data consistency and validity criteria that support the software's robustness.
	\item \emph{Unified OPF formulation and convex relaxation:}
		Based on this system model, we discuss an OPF formulation that supports cost minimization, loss minimization, and a combination of both. The unified model for AC and DC subgrids enables a compact formulation in which the hybrid nature is only explicit in the grid's state, i.e., voltage \emph{phasors} for AC and voltage \emph{magnitudes} for DC subgrids. To simplify the further study of the OPF problem, we introduce a \emph{state space relaxation} that relaxes the state of DC subgrids to voltage \emph{phasors} and prove its \emph{exactness}, i.e., that the optimal DC voltage magnitudes can always be recovered. Therewith, a \emph{unified OPF formulation} is obtained that, in contrast to other unified formulations like~\cite{Bahrami2017a,Baradar2013b} which explicitly retain real-valued state variables in the DC systems (cf.~\cite[Sec.~II-A1]{Bahrami2017a} and~\cite[Eq.~(4)]{Baradar2013b}), additionally unifies the state representation. It thus establishes a \emph{complete analogy} of AC and DC subgrids in the primal and Lagrangian dual domain that enables the \emph{direct} generalization of results on AC grids to DC grids. This is illustrated for the semidefinite and second-order cone relaxation as well as associated results on exactness and locational marginal prices, which can be generalized directly without explicitly considering the hybrid nature of the grid.
	\item \emph{Software design and validation:}
		To arrive at a software framework that is easy to use, comprehend, and adapt, we developed a simple yet extensible object-oriented software design for \hynet, whose fundamental elements are presented here. The simplicity arises from a systematic partitioning of the data flow into object relations, while extensibility is provided via inheritance and standard design patterns. Complementary, we introduced a specially designed relational database schema alongside \hynet to offer an adequate data format for hybrid AC/DC power systems for use in research and education. It is based on the established \emph{structured query language} (SQL) and enables the platform-independent storage of infrastructure and scenario data.
		Further, we illustrate the application of \hynet in several case studies to validate its functionality and showcase its robust and competitive performance.
\end{itemize}
After the submission of this manuscript, another open-source OPF software for hybrid AC/DC grids was independently proposed in~\cite{Ergun2019a}. It focuses on a comprehensive formulation that also supports meshed MT-HVDC systems. In contrast, this work presents a flexible, concise, and unified formulation that promotes ease of use and establishes a beneficial platform for extensions and further mathematical studies.

\subsection{Outline}

In the following, Section~\ref{sec:model} and~\ref{sec:constraints} present the system model for hybrid AC/DC power systems. Section~\ref{sec:opf} discusses the OPF formulation and the state space relaxation. Section~\ref{sec:relaxation} continues with the convex relaxation of the unified OPF problem and discusses the necessitated bus voltage recovery as well as locational marginal prices and exactness of the relaxation. Section~\ref{sec:implementation} highlights \hynet's software design, which is followed by the software validation and case studies in Section~\ref{sec:validation}. Finally, Section~\ref{sec:conclusion} concludes the paper.

\begin{figure*}[!t]%
\centering%
\hfill%
\subfloat[]{%
\includegraphics{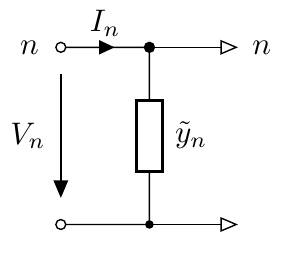}%
\label{fig:bus}}%
\hfill\hfill%
\subfloat[]{%
\includegraphics{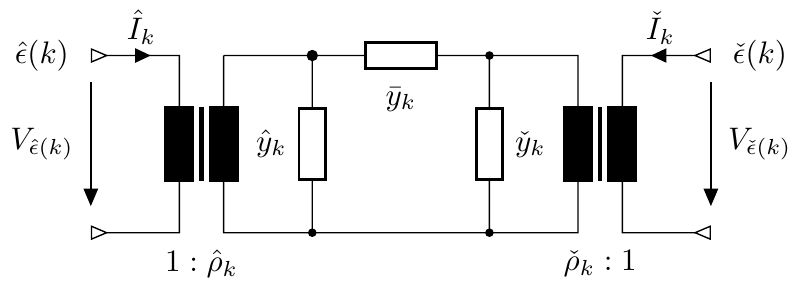}%
\label{fig:branch}}%
\hfill\hfill%
\subfloat[]{%
\includegraphics{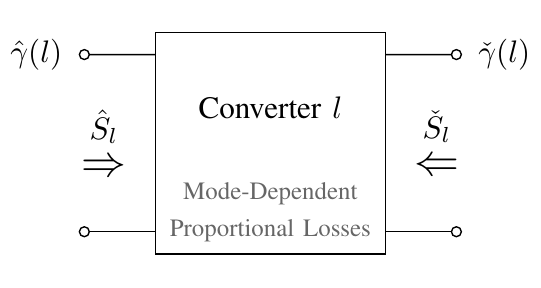}%
\label{fig:converter}}%
\hfill%
\caption{Electrical models: (a) Model for bus $n\in\sV$~\cite{Hotz2016a} with an injection port (circles) and a branch port (triangles). (b) Model for branch $k\in\sE$~\cite{Hotz2016a}, which connects the branch port of source bus $\eSrc(k)$ to the branch port of destination bus $\eDst(k)$. (c) Model for converter $l\in\sC$, which connects the injection port of source bus $\cSrc(l)$ to the injection port of destination bus $\cDst(l)$.}
\end{figure*}

\subsection{Notation}

The set of natural numbers is denoted by $\nN$, the set of integers by $\nZ$, the set of real numbers by $\nR$, the set of nonnegative real numbers by $\nRnn$, the set of positive real numbers by $\nRp$, the set of complex numbers by $\nC$, and the set of Hermitian matrices in $\nC^{N\times N}$ by $\nS^N$. The imaginary unit is denoted by $\displaystyle\iu=\sqrt{-1}$. For $x\in\nC$, its real part is $\real(x)$, its imaginary part is $\imag(x)$, its absolute value is $\abs{x}$, the principal value of its argument is $\arg(x)\in(-\pi,\pi]$, and its complex conjugate is $x\conj$. For a matrix $\mA$, its transpose is $\mA\tran$, its conjugate (Hermitian) transpose is $\mA\herm$, its trace is $\trace(\mA)$, its rank is $\rank(\mA)$, its Frobenius norm is $\norm{\mA}_\mathrm{F}$, and its element in row~$i$ and column~$j$ is $[\mA]_{i,j}$. For two matrices $\mA,\mB\in\nS^N$, $\mA\psd\mB$ denotes that $\mA-\mB$ is positive semidefinite. For real-valued vectors, inequalities are component-wise. The vector $\ve_n$ denotes the $n$th standard basis vector of appropriate dimension. For a countable set $\sC$, its cardinality is $\abs{\sC}$. For a set $\sN\subset\nN$ and vectors or matrices $\vx_n\in\sS$, with $n\in\sN$, $\vx_{\sN}$ denotes the $\abs{\sN}$-tuple $\vx_{\sN} = (\vx_n)_{n\in\sN}$ and $\sS_{\sN}$ the $\abs{\sN}$-fold Cartesian product $\sS_{\sN}={\prod_{n\in\sN} \sS}$.

\section{System Model}
\label{sec:model}

This section presents a system model for hybrid AC/DC grids with P2P- and MT-HVDC systems, which consists of a description of its \emph{network topology} and an \emph{electrical model}. In the following, AC lines, cables, transformers, and phase shifters are referred to as \emph{AC branches}, DC lines and cables as \emph{DC branches}, inverters, rectifiers, VSCs, and back-to-back converters as \emph{converters}, and points of interconnection, generation injection, and load connection are called \emph{buses}.

\subsection{Network Topology}
\label{sec:model:architecture}

The network topology of the hybrid AC/DC power system, which consists of the interconnection of an arbitrary number of AC and DC subgrids via converters, is described by the directed multigraph $\gG = {\{\sV,\sE,\sC,\eSrc,\eDst,\cSrc,\cDst\}}$, where
\begin{enumerate}
	\item $\sV=\{1,\ldots,\cV\}$ is the set of buses,
	\item $\sE=\{1,\ldots,\cE\}$ is the set of branches,
	\item $\sC=\{1,\ldots,\cC\}$ is the set of converters,
	\item $\eSrc:\sE\rightarrow\sV$ maps a branch to its source bus,
	\item $\eDst:\sE\rightarrow\sV$ maps a branch to its destination bus,
	\item $\cSrc:\sC\rightarrow\sV$ maps a converter to its source bus, and
	\item $\cDst:\sC\rightarrow\sV$ maps a converter to its destination bus.
\end{enumerate}
The directionality of branches and converters is not related to the direction of power flow and can be chosen \emph{arbitrarily}. The buses $\sV$ are partitioned into a set $\sVAC$ of \emph{AC buses} and a set $\sVDC$ of \emph{DC buses}, i.e.,
%
% $\sV = \sVAC\cup\sVDC$ and $\sVAC\cap\sVDC = \emptyset$\,.
\begin{align}
	\sV = \sVAC\cup\sVDC
	\qquad\text{and}\qquad
	\sVAC\cap\sVDC = \emptyset\,.
\end{align}
AC and DC buses must not be connected by a branch, i.e., the branches $\sE$ are partitioned into AC and DC branches,
%
% $\sE = \sEAC\cup\sEDC$ and $\sEAC\cap\sEDC = \emptyset$\,.
\begin{align}
	\sE = \sEAC\cup\sEDC
	\qquad\text{and}\qquad
	\sEAC\cap\sEDC = \emptyset\,.
\end{align}
Accordingly, the set $\sEAC$ of \emph{AC branches} and the set $\sEDC$ of \emph{DC branches} is given by
\begin{subequations}
\begin{align}
	\sEAC &= \{k\in\sE: \eSrc(k),\eDst(k)\in\sVAC\}\\
	\sEDC &= \{k\in\sE: \eSrc(k),\eDst(k)\in\sVDC\}\,.
\end{align}
\end{subequations}

Note that the terminal buses of converters are \emph{not} restricted, i.e., the model supports AC/DC and DC/AC as well as AC and DC back-to-back converters. To support the mathematical exposition later on, some terms and expressions are defined.
\begin{definition}
Consider the directed subgraph $\gG' = \{\sV, \sE,\allowbreak \eSrc, \eDst\}$ with all buses and branches. A \emph{connected component}~\cite{Bondy1976a} in the \emph{underlying}~\cite{Bondy1976a} undirected graph of $\gG'$ is called \emph{subgrid}. A subgrid comprising buses in $\sVAC$ is called \emph{AC subgrid}. A subgrid comprising buses in $\sVDC$ is called \emph{DC subgrid}.
\end{definition}
\begin{definition}
The set $\aSrc{\sBE}(n)\subseteq\sE$ and $\aDst{\sBE}(n)\subseteq\sE$ of branches outgoing and incoming at bus $n\in\sV$, respectively,~is
\begin{subequations}
\begin{align}
	\aSrc{\sBE}(n) &= \{k\in\sE: \eSrc(k) = n\}\\
	\aDst{\sBE}(n) &= \{k\in\sE: \eDst(k) = n\}\,.
\end{align}
\end{subequations}
\end{definition}
\begin{definition}
The set $\aSrc{\sBC}(n)\subseteq\sC$ and $\aDst{\sBC}(n)\subseteq\sC$ of converters outgoing and incoming at bus $n\in\sV$, respectively,~is
\begin{subequations}
\begin{align}
	\aSrc{\sBC}(n) &= \{l\in\sC: \cSrc(l) = n\}\\
	\aDst{\sBC}(n) &= \{l\in\sC: \cDst(l) = n\}\,.
\end{align}
\end{subequations}
\end{definition}
Finally, a generally valid property of the network topology is established, which is utilized later on.
\begin{definition}[Self-Loop Free Network Graph]\label{def:noselfloops}
The multigraph $\gG$ does not comprise any self-loops, i.e.,
\begin{align}
	\nexists\,k\in\sE:\eSrc(k)=\eDst(k)
	\quad\text{and}\quad
	\nexists\,l\in\sC:\cSrc(l)=\cDst(l)\,.
\end{align}
\end{definition}

\subsection{Electrical Model}
\label{sec:model:electrical}

The electrical model for branches and buses is adopted from~\cite{Hotz2016a} and extended with DC subgrids as well as a novel converter model. The characterization of generation and load is adopted from~\cite{Hotz2018a} and generalized to the concept of \emph{injectors}.

\subsubsection{Branch Model}

Branches are represented via the common branch model in Fig.~\ref{fig:branch}. For branch $k\in\sE$, it comprises two shunt admittances $\aSrc{y_k},\aDst{y_k}\in\nC$, a series admittance $\aBar{y_k}\in\nC$, and two complex voltage ratios $\aSrc{\rho_k},\aDst{\rho_k}\in\nC\setminus\{0\}$. In the latter, $\abs{\aSrc{\rho_k}}$ and $\abs{\aDst{\rho_k}}$ is the tap ratio and $\arg(\aSrc{\rho_k})$ and $\arg(\aDst{\rho_k})$ the phase shift of the respective transformer, while
\begin{align}
	\rho_k = \aSrc{\rho_k\conj}\aDst{\rho_k}
\end{align}
denotes the \emph{total voltage ratio}. Let the \emph{bus voltage vector} $\vv$, \emph{source current vector} $\viSrc$, and \emph{destination current vector} $\viDst$ be
\begin{align}
	\vv &= [V_1,\ldots,V_\cV]\tran\in\nC^\cV\\
	\viSrc &= [\hskip1.5pt\aSrc{I_1},\ldots,\hskip1pt\aSrc{I_\cE}\hskip1pt]\tran\in\nC^\cE\\
	\viDst &= [\hskip1.5pt\aDst{I_1},\ldots,\hskip1pt\aDst{I_\cE}\hskip1pt]\tran\in\nC^\cE\,.
\end{align}
They are related by Kirchhoff's and Ohm's law, which renders
\begin{align}
	\viSrc = \aSrc{\mY}\vv
	\qquad\text{and}\qquad
	\viDst = \aDst{\mY}\vv
\end{align}
where $\aSrc{\mY},\aDst{\mY}\in\nC^{\cE\times\cV}$ are given in~\cite[Eq.~(6) and~(7)]{Hotz2016a}. In the following, bus voltages are used as \emph{state variables}. While AC subgrids exhibit complex-valued effective (rms) voltage phasors, DC subgrids exhibit real-valued voltages. This is considered by restricting $\vv$ to $\sU$, where
\begin{align}
	\sU &= \{\vv\in\nC^\cV: V_n\in\nRnn,\ n\in\sVDC\}\,.
\end{align}
Furthermore, DC lines and cables are modeled via their series resistance, which is captured as follows.
\begin{definition}\label{def:dc:seriescond}
DC branches equal a series conductance, i.e.,
\begin{align}
	\forall k\in\sEDC&:\,\
		\aSrc{\rho_k}=\aDst{\rho_k}=1,\ \
		\aSrc{y_k}=\aDst{y_k}=0,\ \
		\imag(\aBar{y_k})=0\,.
\end{align}
\end{definition}
Finally, a generally valid physical property of DC branches is observed, which is utilized later on.
\begin{definition}[Lossy DC Branches]\label{def:lossydcbranches}
The series conductance of all DC branches is positive, i.e., $\real(\aBar{y_k}) > 0$, $\forall k\in\sEDC$.
\end{definition}

\subsubsection{Bus Model}

Buses are modeled as depicted in Fig.~\ref{fig:bus}. For bus $n\in\sV$, it comprises a shunt admittance $\aSh{y_n}\in\nC$, connections to the outgoing branches $k\in\aSrc{\sBE}(n)$ as well as to the incoming branches $k\in\aDst{\sBE}(n)$, and an injection port. Let the \emph{injection current vector} $\vi$ be
\begin{align}
	\vi = [I_1,\ldots,I_\cV]\tran\in\nC^\cV\,.
\end{align}
It is related to $\vv$ by Kirchhoff's and Ohm's law, i.e.,
\begin{align}
	\vi = \mY\vv
\end{align}
where the bus admittance matrix $\mY\in\nC^{\cV\times\cV}$ is given in~\cite[Eq.~(10)]{Hotz2016a}. Finally, note that the shunt $\aSh{y_n}$ usually models reactive power compensation and is irrelevant in DC subgrids.
\begin{definition}\label{def:dc:noshunt}
DC buses exhibit a zero shunt admittance, i.e., $\aSh{y_n}=0$, $\forall n\in\sVDC$.
\end{definition}

\subsubsection{Converter Model}

\begin{figure}[!t]%
\vspace{-0.8em}%
\centering%
\footnotesize%
\subfloat[\hspace*{7.5mm}]{%
\includegraphics{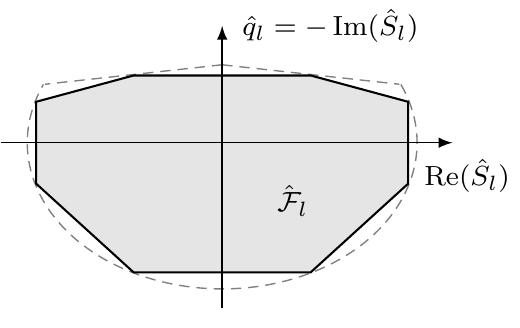}%
\label{fig:pq:converter}}%
\hspace*{0.6mm}%
\subfloat[\hspace*{6mm}]{%
\includegraphics{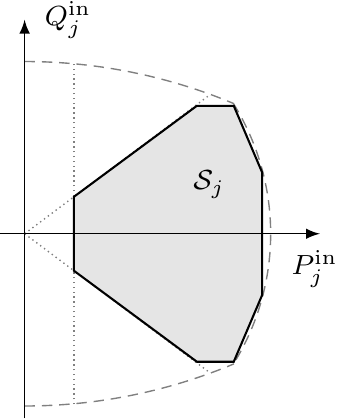}%
\label{fig:pq:generator}\hspace*{-0.6mm}}%
\caption{Qualitative example of the P/Q-capability (dashed) of (a) a typical \emph{voltage source converter}~\cite[Sec.~2.1]{Davies2009a}, \cite[Sec.~2.7.1]{Cigre2014b}, which includes voltage and current related limits, and (b) a typical \emph{generator}~\cite[Ch.~5.4]{Kundur1994a}, which includes physical limits as well as a power factor and minimum output limit (dotted).
The inner polyhedral approximations (solid), which illustrate an exemplary parameterization of the $8$ half-spaces supported in \hynet, restrict the permitted set points to a safe operating area.}%
\label{fig:pq}%
\end{figure}

Converters are modeled as illustrated in Fig.~\ref{fig:converter}. Complementary, the transformer, filter, and phase reactor of a VSC station can be modeled using two AC branches that connect the converter via a filter bus to the point of common coupling, see also~\cite{Zhao2017a,Cigre2014b}. For converter $l\in\sC$, the model considers the source and destination apparent power flow $\aSrc{S_l},\aDst{S_l}\in\nC$, respectively, where active power is converted with a forward and backward conversion loss factor $\aFwd{\eta_l},\aBwd{\eta_l}\in[0,1)$, respectively, while reactive power may be provided, i.e.,
\begin{subequations}\label{eqn:converter:flow}
\begin{align}
	\aSrc{S_l} &= \aSrc{p_l} - (1-\aBwd{\eta_l})\aDst{p_l} - \iu\aSrc{q_l}\\
	\aDst{S_l} &= \aDst{p_l} - (1-\aFwd{\eta_l})\aSrc{p_l} - \iu\aDst{q_l}\,.
\end{align}
\end{subequations}
Therein, $\aSrc{p_l}\in\nRnn$ and $\aDst{p_l}\in\nRnn$ is the nonnegative active power flow from bus $\cSrc(l)$ to $\cDst(l)$ and vice versa, respectively, while $\aSrc{q_l},\aDst{q_l}\in\nC$ is the reactive power support and $\aFwd{\eta_l},\aBwd{\eta_l}$ enable a mode-dependent loss parameterization. The P/Q-capability of the converter at the source and destination bus is approximated by the compact polyhedral set $\aSrc{\sF_l}$ and $\aDst{\sF_l}$, respectively, cf. Fig.~\ref{fig:pq:converter}. In this model, this is captured~by
\begin{align}
	\vf_l = [\aSrc{p_l},\aDst{p_l},\aSrc{q_l},\aDst{q_l}]\tran\in\sF_l\subset\nRnn^2\times\nR^2
\end{align}
where the compact polyhedral set $\sF_l$ is a reformulation of $\aSrc{\sF_l}$ and $\aDst{\sF_l}$ in terms of the \emph{converter state vector} $\vf_l$ using~\eqref{eqn:converter:flow}. For a concise notation, the state vectors of all converters are stacked, i.e.,
\begin{align}
	\vf = [\vf_1\tran,\ldots,\vf_\cC\tran]\tran\in\sF=\prod_{l\in\sC}\sF_l
\end{align}
where the polyhedral set $\sF$ is expressed as
\begin{align}\label{eqn:model:converter:pqcap}
	\sF=\{\vf\in\nR^{4\cC}:\mH\vf\leq\vh\}\,.
\end{align}
In the latter, $\mH\in\nR^{F\times 4\cC}$ and $\vh\in\nR^F$ capture the $F\in\nN$ inequality constraints that describe the capability regions of all converters. Finally, the absence of reactive power on the DC-side of a converter is established.
\begin{definition}\label{def:dc:conv:injection}
The DC-side of all converters exclusively injects active power, i.e.,
\begin{subequations}
\begin{align}
	\forall l\in\bigcup_{n\in\sVDC}\aSrc{\sBC}(n)&:\quad\sF_l\subset\nRnn^2\times\{0\}\times\nR\\
	\forall l\in\bigcup_{n\in\sVDC}\aDst{\sBC}(n)&:\quad\sF_l\subset\nRnn^2\times\nR\times\{0\}\,.
\end{align}
\end{subequations}
\end{definition}

The presented converter model is linear in the converter state vector, rendering it easily tractable. This is attained by excluding the nonlinear mode complementarity constraint
\begin{align}
	\aSrc{p_l}\aDst{p_l} = 0
\end{align}
which can cause a nonnegative and bounded loss error
\begin{align}
	\vartheta_l =
	    \aFwd{\eta_l}\big(\aSrc{p_l} - [\real(\aSrc{S_l})]_{+}\big)
	  + \aBwd{\eta_l}\big(\aDst{p_l} - [\real(\aDst{S_l})]_{+}\big)
\end{align}
in converter $l\in\sC$, where $[x]_{+} = \max(x,0)$. In general, OPF problems do not incentivize an increase of load and, as the loss error (load increase) is minimal (zero) at mode complementarity, this relaxation typically remains without~impact.

In the software framework, the loss error is monitored and, if it exceeds a certain tolerance, the converter mode is fixed according to the net active power flow to ensure a zero loss error and the computation is reissued. Static (no-load) converter losses are modeled explicitly for ease of use, while w.l.o.g. the mathematical model considers them as fixed loads. The P/Q-capabilities $\aSrc{\sF_l}$ and $\aDst{\sF_l}$ are specified by an intersection of up to $8$ half-spaces, offering adequate accuracy at a moderate number of constraints and parametrization~complexity.

\subsubsection{Generators, Prosumers, and Loads}

For a compact notation, it is observed that, in the context of optimal power flow, power producers, prosumers (and flexible distribution systems), as well as flexible and fixed loads are conceptually equivalent. They comprise a certain set of valid operating points in the P/Q-plane and quantify their preferences via a real-valued function over the P/Q-plane. This is utilized for an abstraction of these entities to injectors. An \emph{injector} can inject a certain amount of (positive or negative) active and reactive power, which is associated with a certain cost. The set of injectors is denoted by $\sI=\{1,\ldots,\cI\}$. Injector $j\in\sI$ injects the active power $\aIj{P_j}$ and the reactive power $\aIj{Q_j}$, which are collected in the \emph{injection vector} $\vs_j=[\aIj{P_j}, \aIj{Q_j}]\tran\in\nR^2$. Its valid operating points are specified by the \emph{capability region} $\sS_j\subset\nR^2$, which is a nonempty, compact, and convex set. The convex \emph{cost function} $\fC_j:\sS_j\rightarrow\nR$ specifies the cost associated with an operating point. The injector's terminal bus is specified by $\nSrc:\sI\rightarrow\sV$, i.e., injector $j\in\sI$ is connected to the injection port of bus $\nSrc(j)$. Conversely, the injectors at a bus are identified as follows.
\begin{definition}
The set $\aSrc{\sBI}(n)\subseteq\sI$ of injectors connected to bus $n\in\sV$ is
\begin{align}
	\aSrc{\sBI}(n) = \{j\in\sI: \nSrc(j) = n\}\,.
\end{align}
\end{definition}
Moreover, the nature of DC subgrids is respected.
\begin{definition}\label{def:dc:injector:injection}
Injectors connected to DC buses exclusively inject active power, i.e., $\sS_j\subset\nR\times\{0\}$, $\forall j\in\bigcup_{n\in\sVDC}\aSrc{\sBI}(n)$.
\end{definition}

For example, for a generator, $\sS_j$ is a convex approximation of its P/Q-capability, while $\fC_j$ reflects the generation cost. For a fixed load, $\sS_j$ is singleton and $\fC_j$ maps to a constant value. For a flexible load, $\sS_j$ characterizes the implementable load shift, while $\fC_j$ reflects the cost for load dispatching.

In the software framework, fixed loads are modeled explicitly for ease of use. Furthermore, the P/Q-capability $\sS_j$ is specified as a polyhedral set defined by the intersection of up to $8$ half-spaces, which can approximate the physical capabilities and restrict the power factor, cf. Fig.~\ref{fig:pq:generator}. The cost function $\fC_j:\sS_j\subset\nR^2\rightarrow\nR$ is considered linearly separable in the active and reactive power costs, i.e.,
\begin{align}
	\fC_j(\vs) = \aP{\fC_j}(\ve_1\tran\vs) + \aQ{\fC_j}(\ve_2\tran\vs)\,.
\end{align}
The active and reactive power cost functions $\aP{\fC_j},\aQ{\fC_j}:\nR\rightarrow\nR$ are considered convex and \emph{piecewise linear}. Hence, they can approximate any convex function with arbitrary accuracy.

\subsubsection{Power Balance}

The flow conservation arising from Kirchhoff's current law balances the nodal injections with the flow into branches and converters. This is captured by the \emph{power balance} equations
\begin{subequations}\label{eqn:model:crt:powbal}
\begin{align}
	\vv\herm\mP_n\vv + \vp_n\tran\vf &= \ve_1\tran\sum_{\mathclap{j\in\aSrc{\sBI}(n)}} \vs_j\,,
		&& \forall n\in\sV
		\label{eqn:model:crt:powbal:p}\\[0.2em]
	\vv\herm\mQ_n\vv + \vq_n\tran\vf &= \ve_2\tran\sum_{\mathclap{j\in\aSrc{\sBI}(n)}} \vs_j\,,
		&& \forall n\in\sV\,.
		\label{eqn:model:crt:powbal:q}
\end{align}
\end{subequations}
Therein, the left hand side describes the flow of active and reactive power into the branches and converters, respectively, while the right hand side accumulates the nodal active and reactive power injection. The matrices $\mP_n,\mQ_n\in\nS^\cV$ are a function of the bus admittance matrix and given in~\cite[Eq.~(14)]{Hotz2016a}. The vectors $\vp_n,\vq_n\in\nR^{4\cC}$, which include~\eqref{eqn:converter:flow} and characterize the flow into converters, can be derived as
\begin{align}
	\vp_n &=
	\sum_{l\in\aSrc{\sBC}(n)}\!\!\!\left(\ve_{4l-3}-(1-\aBwd{\eta_l})\ve_{4l-2}\right)\notag\\
	      &\qquad+\!\!
	\sum_{l\in\aDst{\sBC}(n)}\!\!\!\left(\ve_{4l-2}-(1-\aFwd{\eta_l})\ve_{4l-3}\right)\\
	\vq_n &= -\!\!\!\sum_{l\in\aSrc{\sBC}(n)}\!\!\!\ve_{4l-1}\ 
			 -\!\!\!\sum_{l\in\aDst{\sBC}(n)}\!\!\!\ve_{4l}\,. %{4(l-1)+4}
\end{align}

\subsubsection{Electrical Losses}
\label{sec:model:electrical:losses}

The total electrical losses amount to the difference of total active power generation and load, i.e., the sum of the right-hand side in~\eqref{eqn:model:crt:powbal:p} over all $n\in\sV$ or, equivalently, the respective sum of the left-hand side. With the latter and $\vv\herm\mA\vv=\trace(\mA\vv\vv\herm)$, the \emph{total electrical losses} can be derived as $\fLoss(\vv\vv\herm,\vf)$, where $\fLoss:\nS^\cV\times\nR^{4\cC}\rightarrow\nR$ is\footnote{The motivation for defining $\fLoss$ in terms of the outer product of the bus voltage vector will become evident in Section~\ref{sec:relaxation}.}
\begin{align}
	\fLoss(\mV,\vf) = \trace(\mL\mV) + \vl\tran\vf
\end{align}
with $\mL\in\nS^\cV$ and $\vl\in\nR^{4\cC}$ given by
\begin{align}\label{eqn:model:loss:coeff}
	\mL = \frac{1}{2}(\mY + \mY\herm)
	\quad\text{and}\quad
	\vl = \sum_{l\in\sC} \left[\aFwd{\eta_l}\ve_{4l-3} + \aBwd{\eta_l}\ve_{4l-2}\right].
\end{align}

\section{System Constraints}
\label{sec:constraints}

In the following, the formulation of physical and stability-related limits is presented. It is based on~\cite{Hotz2016a}, where the common apparent power flow limit (``MVA rating'') is substituted by its underlying ampacity, voltage drop, and angle difference constraint (cf.~\cite[Ch.~6.1.12]{Kundur1994a}, \cite[Ch.~4.9]{Bergen2000a}) to improve expressiveness and mathematical structure.

\setcounter{subsubsection}{0}

\subsubsection{Voltage}

Due to physical and operational requirements, the voltage at bus $n\in\sV$ must satisfy $\abs{V_n}\in[\aLB{V_n},\aUB{V_n}]\subset\nRnn$, where $\aUB{V_n}>\aLB{V_n}>0$. With $\mM_n = \ve_n\ve_n\tran\in\nS^\cV$, this reads
\begin{align}
	\aLB{V_n}^2 \leq \vv\herm\mM_n\vv \leq \aUB{V_n}^2\,,
	\qquad \forall n\in\sV\,.
		\label{eqn:crt:voltage}
\end{align}

\subsubsection{Ampacity}

The thermal flow limit on branch $k\in\sE$ can be expressed as $\abs{\aSrc{I_k}} \leq \aSrcUB{I}_k$ and $\abs{\aDst{I_k}} \leq \aDstUB{I}_k$, where $\aSrcUB{I}_k,\aDstUB{I}_k\in\nRp$. In quadratic form, this renders
\begin{align}
	\vv\herm\mISrc_k\vv \leq \aSrcUB{I}_k^2\,,
	&&
	\vv\herm\mIDst_k\vv \leq \aDstUB{I}_k^2\,,
	&& \forall k\in\sE
		\label{eqn:crt:ampacity}
\end{align}
where $\mISrc_k,\mIDst_k\in\nS^\cV$ are given in~\cite[Eq.~(18)]{Hotz2016a}.

\subsubsection{Voltage Drop}

The stability-related limit on the voltage drop $\nu_k\in\nR$ along AC branch $k\in\sEAC$, i.e.,
\begin{equation}
	\nu_k = \abs{V_{\eDst(k)}}/\abs{V_{\eSrc(k)}} - 1
\end{equation}
is $\nu_k\in[\aLB{\nu_k},\aUB{\nu_k}]\subset[-1,\infty)$, with $\aLB{\nu_k}<\aUB{\nu_k}$. This is captured by
\begin{align}
	\vv\herm\aLB{\mM_k}\vv \leq 0\,,
	&&
	\vv\herm\aUB{\mM_k}\vv \leq 0\,,
	&& \forall k\in\sEAC
		\label{eqn:crt:drop}
\end{align}
in which $\aLB{\mM_k},\aUB{\mM_k}\in\nS^\cV$ are given in~\cite[Eq.~(24) and~(25)]{Hotz2016a}.

\subsubsection{Angle Difference}

The stability-related limit on the voltage angle difference $\delta_k\in\nR$ along AC branch $k\in\sEAC$, i.e.,
\begin{equation}
	\delta_k = \arg(V_{\eSrc(k)}\conj V_{\eDst(k)})
\end{equation}
reads $\delta_k\in[\aLB{\delta_k},\aUB{\delta_k}]\subset(-\pi/2,\pi/2)$, with $\aLB{\delta_k}<\aUB{\delta_k}$. Note that
\begin{align}
	\real(V_{\eSrc(k)}\conj &V_{\eDst(k)}) \geq 0
		\hspace*{-0.2em}
		\label{eqn:crt:angle:realpart}\\
	\tan(\aLB{\delta_k}) \leq \imag(V_{\eSrc(k)}\conj V_{\eDst(k)})\,/\,&\real(V_{\eSrc(k)}\conj V_{\eDst(k)}) \leq \tan(\aUB{\delta_k})
		\label{eqn:crt:angle:tanpart}
\end{align}
is an equivalent formulation of this constraint that can be put~as
\begin{align}
	\vv\herm\mA_k\vv \leq 0\,,
		\ \
	\vv\herm\aLB{\mA_k}\vv &\leq 0\,,
		\ \
	\vv\herm\mAUB_k\vv \leq 0\,,
		\ \ \forall k\in\sEAC
		\label{eqn:crt:angle}
\shortintertext{where}
	\mA_k &= -\aSrc{\mM_k} - \aSrc{\mM_k}\herm
		\label{eqn:crt:mtx:ak}
\end{align}
with $\aSrc{\mM_k} = \ve_{\eSrc(k)}\ve_{\eDst(k)}\tran\in\nR^{\cV\times\cV}$ and $\aLB{\mA_k},\mAUB_k\in\nS^\cV$ as given in~\cite[Eq.~(30) and~(31)]{Hotz2016a}.

\section{Optimal Power Flow}
\label{sec:opf}

The OPF problem identifies the optimal utilization of the grid infrastructure and generation resources to satisfy the load, where optimality is typically considered with respect to minimum injection costs or minimum electrical losses. With the system model and system constraints above, the OPF problem can be cast as the following optimization problem.%
\begin{subequations}\label{eqn:opf:nonconvex:full}%
\begin{alignat}{2}
&\minimize_{
	\mathclap{\substack{\\[0.1em]
		\vs_j\in\sS_j,\,\vv\in\sU\\[0.1em]
		\vf\in\sF}}}
	\quad\ &&
	\sum_{j\in\sI} \fC_j(\vs_j) + \tau\fLoss(\vv\vv\herm,\vf)\\[0.25em]
&\subjectto &&
\eqref{eqn:model:crt:powbal},
\eqref{eqn:crt:voltage},
\eqref{eqn:crt:ampacity},
\eqref{eqn:crt:drop},
\eqref{eqn:crt:angle}\,.\
	\label{eqn:opf:nonconvex:full:crt}
\end{alignat}
\end{subequations}

The objective consists of the injection costs and a penalty term comprising the electrical losses weighted by an (artificial) \emph{loss price} $\tau\geq 0$, which enables injection cost minimization, electrical loss minimization, and a combination of both. The penalty term offers a convenient parameterization, while the same objective can also be expressed with cost functions only.
\begin{proposition}
Consider the OPF problem~\eqref{eqn:opf:nonconvex:full}. Including a loss penalty with $\tau>0$ is equivalent to increasing the marginal cost of active power by $\tau$ for all injectors $j\in\sI$.
\end{proposition}
\begin{IEEEproof}
Considering Section~\ref{sec:model:electrical:losses}, it follows that
\begin{equation*}
\begin{IEEEeqnarraybox}{rCl}
	\sum_{j\in\sI} \fC_j(\vs_j) + \tau\fLoss(\vv\vv\herm,\vf)
		&=& \sum_{j\in\sI} \fC_j(\vs_j) + \tau\sum_{n\in\sV} \ve_1\tran\sum_{\mathclap{j\in\aSrc{\sBI}(n)}} \vs_j\notag\\[0.1em]
		&=& \sum_{j\in\sI} \left[ \fC_j(\vs_j) + \tau\ve_1\tran\vs_j \right]\,.
		\vspace{-1.4em}
\end{IEEEeqnarraybox}
\vspace{0.25em}
\IEEEQEDhereeqn
\end{equation*}
\end{IEEEproof}

In~\eqref{eqn:opf:nonconvex:full}, it can be observed that the objective and constraints consider AC and DC subgrids in a unified manner, while their bus voltages are treated differently by restricting $\vv$ to $\sU$. With respect to the Lagrangian dual domain and convex relaxations, this restriction complicates further mathematical studies. To avoid these issues, it is observed that all currently installed and almost all planned HVDC systems are P2P-HVDC or \emph{radial} MT-HVDC systems~\cite{Buigues2017a}. This observation is formalized by the following definition and, in the next section, it is utilized to unify the representation of AC and DC voltages.
\begin{requirement}[Radial DC Subgrids]\label{def:dcgrid:radial}
The underlying undirected graph of the directed subgraph $\aDC{\gG} = {\{\sVDC, \sEDC, \eSrc, \eDst\}}$ is \emph{acyclic}~\cite{Bondy1976a}, i.e., its connected components are~\emph{trees}~\cite{Bondy1976a}.
\end{requirement}

\subsection{State Space Relaxation}
\label{sec:opf:unified}

To eliminate the restriction of $\vv$ to $\sU$ and, therewith, \emph{unify} the representation of AC and DC voltages in the OPF problem, let~\eqref{eqn:crt:angle:realpart} also be imposed on DC subgrids, i.e.,
\begin{align}\label{eqn:crt:angle:dcb}
	\vv\herm\mA_k\vv \leq 0\,,
		\qquad\forall k\in\sEDC
\end{align}
with $\mA_k\in\nS^\cV$ in~\eqref{eqn:crt:mtx:ak}. Jointly, Def.~\ref{def:dcgrid:radial} and the complementary constraints in~\eqref{eqn:crt:angle:dcb} enable the following result.
\begin{theorem}\label{thm:dcvoltages}
Consider any $\vf\in\sF$ and $\vs_j\in\sS_j$, for $j\in\sI$. Let $\vv\in\nC^\cV$ satisfy the constraints in~\eqref{eqn:opf:nonconvex:full:crt} and~\eqref{eqn:crt:angle:dcb} for the given $\vf$ and $\vs_\sI$. Then, the bus voltage vector $\aUBar{\vv}\in\sU$ given by
\begin{align}
	[\aUBar{\vv}]_n =
		\begin{cases}
			\hphantom{|}[\vv]_n & \text{if}\ n\in\sVAC\\
			\abs{[\vv]_n}       & \text{if}\ n\in\sVDC
		\end{cases}
	\label{eqn:dcvoltages:reconst}
\end{align}
satisfies~\eqref{eqn:opf:nonconvex:full:crt} and~\eqref{eqn:crt:angle:dcb} for the given $\vf$ and $\vs_\sI$ with equivalent constraint function values and $\fLoss(\aUBar{\vv}\aUBar{\vv}\herm,\vf)=\fLoss(\vv\vv\herm,\vf)$.
\end{theorem}
\begin{IEEEproof}
See Appendix~\ref{apx:thm:dcvoltages}.
\end{IEEEproof}

Considering that $\sU\subset\nC^\cV$, Theorem~\ref{thm:dcvoltages} enables an \emph{exact state space relaxation} of the OPF problem~\eqref{eqn:opf:nonconvex:full} to $\vv\in\nC^\cV$ under Def.~\ref{def:dcgrid:radial} and the complementary constraints~\eqref{eqn:crt:angle:dcb}. Thus, for radial DC subgrids the representation of AC and DC voltages can be \emph{unified}, where the corresponding OPF formulation~reads%
\begin{subequations}\label{eqn:opf:nonconvex:unified:full}%
\begin{alignat}{2}
\aOpt{p} =\;
&\minimize_{
	\mathclap{\substack{\\[0.1em]
		\vs_j\in\sS_j,\,\vv\in\nC^\cV\\[0.1em]
		\vf\in\sF}}}
	\quad\ &&
	\sum_{j\in\sI} \fC_j(\vs_j) + \tau\fLoss(\vv\vv\herm,\vf)\\[0.25em]
&\subjectto &&
\eqref{eqn:model:crt:powbal},
\eqref{eqn:crt:voltage},
\eqref{eqn:crt:ampacity},
\eqref{eqn:crt:drop},
\eqref{eqn:crt:angle},
\eqref{eqn:crt:angle:dcb}\,.
	\label{eqn:opf:nonconvex:unified:full:crt}
\end{alignat}
\end{subequations}

\subsection{Unified Optimal Power Flow Formulation}

Concluding, for a concise statement of the unified OPF formulation the constraints in~\eqref{eqn:opf:nonconvex:unified:full:crt} and the restriction of $\vf$ to $\sF$ are expressed as ${(\vv\vv\herm,\vf,\vs_\sI)\in\sX}$, where
\begin{subequations}\label{eqn:opf:nonconvex:unified:compact:xdef}
\begin{align}
	\sX = \Big\{\,
		&(\mV,\vf,\vs_\sI) \in \nS^\cV\times\nR^{4\cC}\times\nR^2_\sI :\hspace{-0.75em}
		\\[0.25em]
		&\trace(\mP_n\mV) + \vp_n\tran\vf = \ve_1\tran\sum_{\mathclap{j\in\aSrc{\sBI}(n)}} \vs_j\,,
			&& \forall n\in\sV
			\label{eqn:opf:nonconvex:unified:compact:powbal:p}\\[0.2em]
		&\trace(\mQ_n\mV) + \vq_n\tran\vf = \ve_2\tran\sum_{\mathclap{j\in\aSrc{\sBI}(n)}} \vs_j\,,
			&& \forall n\in\sV
			\label{eqn:opf:nonconvex:unified:compact:powbal:q}\\
		&\trace(\mC_m\mV) + \vc_m\tran\vf \leq b_m\,,
			&& \forall m\in\sM
			\label{eqn:opf:nonconvex:unified:compact:crt:ieq}
	\,\Big\}\,.
\end{align}
\end{subequations}
Therein, \eqref{eqn:opf:nonconvex:unified:compact:crt:ieq} captures \eqref{eqn:model:converter:pqcap}, \eqref{eqn:crt:voltage}, \eqref{eqn:crt:ampacity}, \eqref{eqn:crt:drop}, \eqref{eqn:crt:angle}, and \eqref{eqn:crt:angle:dcb}, where $\mC_m$, $\vc_m$, and $b_m$, $m\in\sM=\{1,\ldots,\cM\}$, reproduce the ${\cM = F + 2\cV + 3\cE + 4\cEAC}$ inequality constraints. With the polyhedral set $\sX$ in~\eqref{eqn:opf:nonconvex:unified:compact:xdef}, \hynet's unified OPF formulation is given by
\begin{subequations}\label{eqn:opf:nonconvex:unified:compact}%
\begin{alignat}{2}
\aOpt{p} =\;
&\minimize_{
	\mathclap{\substack{\\[0.1em]
		\vs_j\in\sS_j,\,\vv\in\nC^\cV\\[0.1em]
		\vf\in\nR^{4\cC}}}}
	\quad\ &&
	\sum_{j\in\sI} \fC_j(\vs_j) + \tau\fLoss(\vv\vv\herm,\vf)
		\label{eqn:opf:nonconvex:unified:compact:objective}\\
&\subjectto &&
	(\vv\vv\herm,\vf,\vs_\sI)\in\sX\,.
\end{alignat}
\end{subequations}

In the software framework, the injector cost functions are included in \emph{epigraph form}~\cite{Boyd2004a} to obtain an explicit formulation as a quadratically-constrained quadratic problem.

\section{Convex Relaxation}
\label{sec:relaxation}

In the following, the semidefinite and second-order cone relaxation of the OPF problem is illustrated. By virtue of the \emph{unified} OPF formulation, they can be applied without explicitly considering the hybrid nature of the grid, enabling a particularly simple and concise derivation and result. Subsequently, the necessitated bus voltage recovery is discussed and some aspects related to locational marginal prices and exactness of the relaxations are highlighted.

\subsection{Semidefinite Relaxation}
\label{sec:relaxation:sdr}

By definition, the objective in~\eqref{eqn:opf:nonconvex:unified:compact} and $\sS_j$ are convex, while $\sX$ is a polyhedral set and, thus, also convex. Consequently, the nonconvexity of~\eqref{eqn:opf:nonconvex:unified:compact} arises from the quadratic dependence on the bus voltage vector $\vv$. The outer product $\vv\vv\herm$ may be expressed equivalently by a Hermitian matrix $\mV\in\nS^\cV$, which is positive semidefinite (psd) and of rank~1. As the set of psd matrices is a convex cone, the nonconvexity is then lumped to the rank constraint. In \emph{semidefinite relaxation} (SDR), the rank constraint is omitted to obtain a convex optimization problem, i.e., the SDR of the OPF problem~\eqref{eqn:opf:nonconvex:unified:compact} reads
\begin{subequations}\label{eqn:opf:sdp}%
\begin{alignat}{2}
\aOpt{\hat{p}} =\;
&\minimize_{
	\mathclap{\substack{\\[0.1em]
		\vs_j\in\sS_j,\,\mV\in\nS^\cV\\[0.1em]
		\vf\in\nR^{4\cC}}}}
	\quad\ &&
	\sum_{j\in\sI} \fC_j(\vs_j) + \tau\fLoss(\mV,\vf)\\
&\subjectto
	&&(\mV,\vf,\vs_\sI)\in\sX\\
	&&&\mV\psd\mZero\,.
		\label{eqn:opf:sdp:vpsd}
\end{alignat}
\end{subequations}

Therewith, the OPF problem gains access to the powerful theory of convex analysis as well as solution algorithms with polynomial-time convergence to a globally optimal solution. The SDR provides a solution to~\eqref{eqn:opf:nonconvex:unified:compact} if the optimizer $\aOpt{\mV}$ obtained from~\eqref{eqn:opf:sdp} has rank~1, in which case the relaxation is called \emph{exact}.

An issue of SDR is the quadratic increase in dimensionality by $\cV^2 - 2\cV$ (real-valued) variables, which impedes computational tractability for large-scale grids. Fortunately, the matrices $\mP_n$, $\mQ_n$, and $\mC_m$ in~\eqref{eqn:opf:nonconvex:unified:compact:xdef} and $\mL$ in~\eqref{eqn:model:loss:coeff} are sparse and may only exhibit a nonzero element in row~$i$ and column~$j$ if $(i,j)\in\sJ$, where
\begin{align}\label{eqn:jdef}
	\sJ = {\{ (i,j)\in\sV\times\sV : i,j\in\{\eSrc(k),\eDst(k)\}, k\in\sE \}}\,,
\end{align}
see~\cite[Appendix~B]{Hotz2016a}. This sparsity can be utilized to effectively solve~\eqref{eqn:opf:sdp} also for large-scale systems, e.g., via a \emph{chordal conversion}~\cite{Fukuda2001a}, see also \cite{Jabr2012a,Molzahn2013a,Andersen2014a,Eltved2019a}.

\begin{figure*}[!t]
\centering
\includegraphics{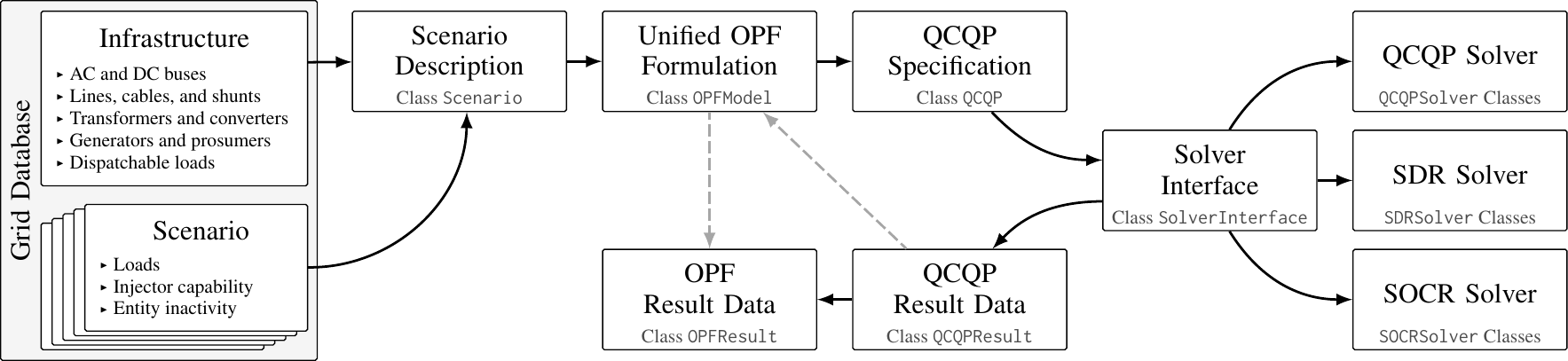}%
\caption{Illustration of \hynet's fundamental design and data flow.}
\label{fig:framework}
\end{figure*}

\subsection{Second-Order Cone Relaxation}
\label{sec:relaxation:socr}

To mitigate the dimensionality uplift of SDR and improve computational efficiency, an additional \emph{second-order cone relaxation} (SOCR) may be applied, see e.g.~\cite{Low2014a,Taylor2015a,Hotz2018a} and the references therein. In the SOCR, the psd constraint~\eqref{eqn:opf:sdp:vpsd} is relaxed to psd constraints on {$2$$\times$$2$} principal submatrices and, as a consequence, only those elements of $\mV$ that are involved in~\eqref{eqn:opf:nonconvex:unified:compact:xdef} and the objective function need to be retained. Correspondingly, the SOCR of~\eqref{eqn:opf:sdp} can be stated as
\begin{subequations}\label{eqn:opf:socp}%
\begin{alignat}{2}
\aOpt{\tilde{p}} =\;
&\minimize_{
	\mathclap{\substack{\\[0.1em]
		\vs_j\in\sS_j,\,\mVp\in\nSp^\cV\\[0.1em]
		\vf\in\nR^{4\cC}}}}
	\quad\ &&
	\sum_{j\in\sI} \fC_j(\vs_j) + \tau\fLoss(\mVp,\vf)\\
&\subjectto
	&&(\mVp,\vf,\vs_\sI)\in\sX\\
	&&&\mS_k\tran\mVp\mS_k\psd\mZero,
		\quad k\in\sE
		\label{eqn:opf:socp:psd}
\end{alignat}
\end{subequations}
in which $\mS_k=[\ve_{\eSrc(k)}, \ve_{\eDst(k)}]\in\nR^{\cV\times 2}$ and $\nSp^\cV$ is the set of Hermitian partial matrices on the graph $\gG' = {\{\sV, \sE, \eSrc, \eDst\}}$, i.e.,
\begin{align}
	\nSp^\cV = \{\fProj(\mV) : \mV\in\nS^\cV\} \subseteq \nS^\cV
\end{align}
with the projection $\fProj$ onto the sparsity pattern given by
\begin{align}\label{eqn:projsparsitypattern}
	\fProj(\mV) = \sum_{(i,j)\in\sJ} [\mV]_{i,j}\ve_i\ve_j\tran\,.
\end{align}
In the SOCR, \eqref{eqn:opf:socp:psd} can be implemented as second-order cone constraints, cf. e.g.~\cite{Kim2003a}. By virtue of the partial matrix, \eqref{eqn:opf:socp} introduces only $2\cE-\cV$ additional (real-valued) variables compared to~\eqref{eqn:opf:nonconvex:unified:compact}.
The SOCR is \emph{exact} if the optimizer $\aOpt{\mVp}$ obtained from~\eqref{eqn:opf:socp} permits a rank-1 completion.

\subsection{Bus Voltage Recovery}
\label{sec:relaxation:recovery}

Considering the potential inexactness of a relaxation and the finite precision of solvers, the bus voltages $\aOpt{\hat{\vv}}$ associated with an optimizer $\aOpt{\mV}$ of the SDR or SOCR must be recovered via a \emph{rank-1 approximation} of $\aOpt{\mV}$ on its \emph{sparsity pattern}.\footnote{Due to voltage-decoupled subgrids, an optimizer $\aOpt{\mV}$ of the SDR~\eqref{eqn:opf:sdp} still allows exact recovery if it permits a decomposition $\aOpt{\mV} = \sum_{i=1}^N \vv_i\vv_i\herm$, where $N$ is the number of subgrids and $\vv_i$ comprises only nonzero elements for bus voltages of subgrid~$i$.} The software framework supports the integration and utilization of different approximation methods. At this point, it includes the method in~\cite[Sec.~III-B-3]{Bose2012b}, which is a computationally efficient graph traversal based approximation, and an approximation in the least-squares sense, i.e.,
\begin{align}\label{eqn:recovery:problem}
	\aOpt{\hat{\vv}} = \arg\min_{\mathclap{\hspace{0.4em}\vv\in\nC^\cV}}\quad
		\norm{\fProj(\vv\vv\herm - \aOpt{\mV})}_\mathrm{F}^2\,.
\end{align}
The (nonconvex) problem~\eqref{eqn:recovery:problem} is solved with a Wirtinger calculus based gradient descent method~\cite{Li2008a} and Armijo's rule for step size control~\cite[Sec.~8.3]{Bazaraa2006a} using an initial point obtained with the method in~\cite[Sec.~III-B-3]{Bose2012b}. 

For a solution $(\aOpt{\mV},\aOpt{\vf},\aOpt{\vs_\sI})$ of the SDR or SOCR with the corresponding rank-1 approximation $\aOpt{\hat{\vv}}$, the software reports the \emph{reconstruction error}
\begin{align}\label{eqn:recovery:error}
	\hat{\kappa}(\aOpt{\mV}) = \norm{\fProj(\aOpt{\hat{\vv}}(\aOpt{\hat{\vv}})\herm - \aOpt{\mV})}_\mathrm{F}^2\,/\,\cJ
\end{align}
as a measure for the (numerical or inherent) inexactness of the relaxation. Further, the induced nodal \emph{power balance error}
\begin{align}
	\hspace{-0.75em}
	\epsilon_n =
		(\aOpt{\hat{\vv}})\herm(\mP_n + \iu\mQ_n)\aOpt{\hat{\vv}}
		+ (\vp_n + \iu\vq_n)\tran\aOpt{\vf}
		- \sum_{\mathclap{j\in\aSrc{\sBI}(n)}} \vw\tran\aOpt{\vs_j}
\end{align}
where $\vw=\left[\,1,\;\iu\,\right]\tran\in\nC^2$, is provided for all $n\in\sV$, which quantifies the impact of the rank-1 approximation on the power flow accuracy (see also the discussion in~\cite[Sec.~III-C]{Molzahn2013a}).

\subsection{Locational Marginal Prices}
\label{sec:relaxation:lmp}

In electricity markets, \emph{nodal pricing} can be applied to account for system constraints and losses and, under \emph{perfect competition} (price takers), the optimal nodal prices are the \emph{locational marginal prices} (LMPs)~\cite{Kirschen2004a}. The LMP at a bus is the cost of serving an increment of load by the cheapest possible injection~\cite{Kirschen2004a}. Thus, it quantifies the sensitivity of the optimal objective value of a cost-minimizing OPF problem with respect to a perturbation of the nodal power balance~\cite{Hotz2018a}. Due to the nonconvexity of the OPF problem, accurate LMPs are in general hard to obtain~\cite{Overbye2004a,Wang2007a}, but in an \emph{exact} SDR and SOCR the LMPs for active and reactive power are given by the optimal Lagrangian dual variables of the power balance constraints~\cite{Hotz2018a}. By virtue of the \emph{unified} OPF formulation, the result in~\cite{Hotz2018a} generalizes\footnote{The proof follows along the lines of~\cite[Sec.~VI-B]{Hotz2018a}, while considering the change to minimization in the primal domain.} to the SDR in~\eqref{eqn:opf:sdp} and the SOCR in~\eqref{eqn:opf:socp}. In \hynet, primal-dual interior-point solvers are applied to the SDR and SOCR. Thus, in case of exactness, the LMPs are obtained as a byproduct of the OPF computation.

\subsection{Exactness of the Relaxations}
\label{sec:relaxation:exactness}

While exactness of a relaxation can be verified \emph{a posteriori} via the reconstruction error in~\eqref{eqn:recovery:error}, \emph{a priori} information about the applicability of a relaxation is desired, i.e., on a system's tendency toward exactness. In the literature, several results were presented for AC grids~\cite{Lavaei2012a,Zhang2013a,Farivar2013a,Low2014b,Madani2015a,Nick2018a}, DC grids~\cite{Gan2014b,Low2014b,Tan2015a}, and hybrid AC/DC grids with P2P-HVDC systems~\cite{Hotz2016a,Hotz2018a}. By virtue of the \emph{unified} OPF formulation, the result in~\cite[Sec.~VII]{Hotz2018a} generalizes\footnote{The proof follows along the lines of~\cite[Sec.~VII]{Hotz2018a}, while considering that the sufficient condition for existence of a psd rank-1 completion in~\cite[Eq.~(25)]{Hotz2018a} generalizes to voltage-decoupled radial subgrids. The latter is evident if the rank-1 completion method in~\cite[Sec.~III-B-3]{Bose2012b} is considered.} under Def.~\ref{def:architecture} to the SOCR~\eqref{eqn:opf:socp} and, thus, to the SDR~\eqref{eqn:opf:sdp}.
\begin{requirement}[Hybrid Architecture]\label{def:architecture}
The underlying undirected graph of the directed subgraph $\gG' = {\{\sV, \sE, \eSrc, \eDst\}}$ is \emph{acyclic}, i.e., its connected components are~\emph{trees}.
\end{requirement}

The result states that for the \emph{hybrid architecture} the SDR and SOCR may only be inexact if the LMPs form a \emph{pathological price profile}, i.e., if they combine to a point in a union of linear subspaces -- a set of measure zero~\cite[Sec.~\mbox{VII-B}]{Hotz2018a}. This suggests that inexactness is unlikely and, if the SDR or SOCR in \hynet is applied to a system with \emph{radial (acyclic) subgrids} like radial distribution grids, exactness may be expected under normal operating conditions.

\begin{table*}[!t]
\caption{Optimal Power Flow Results for the PGLib OPF Benchmark Library for Typical Operating Conditions}
\label{tab:res:pglib:typ}
\centering
\includegraphics{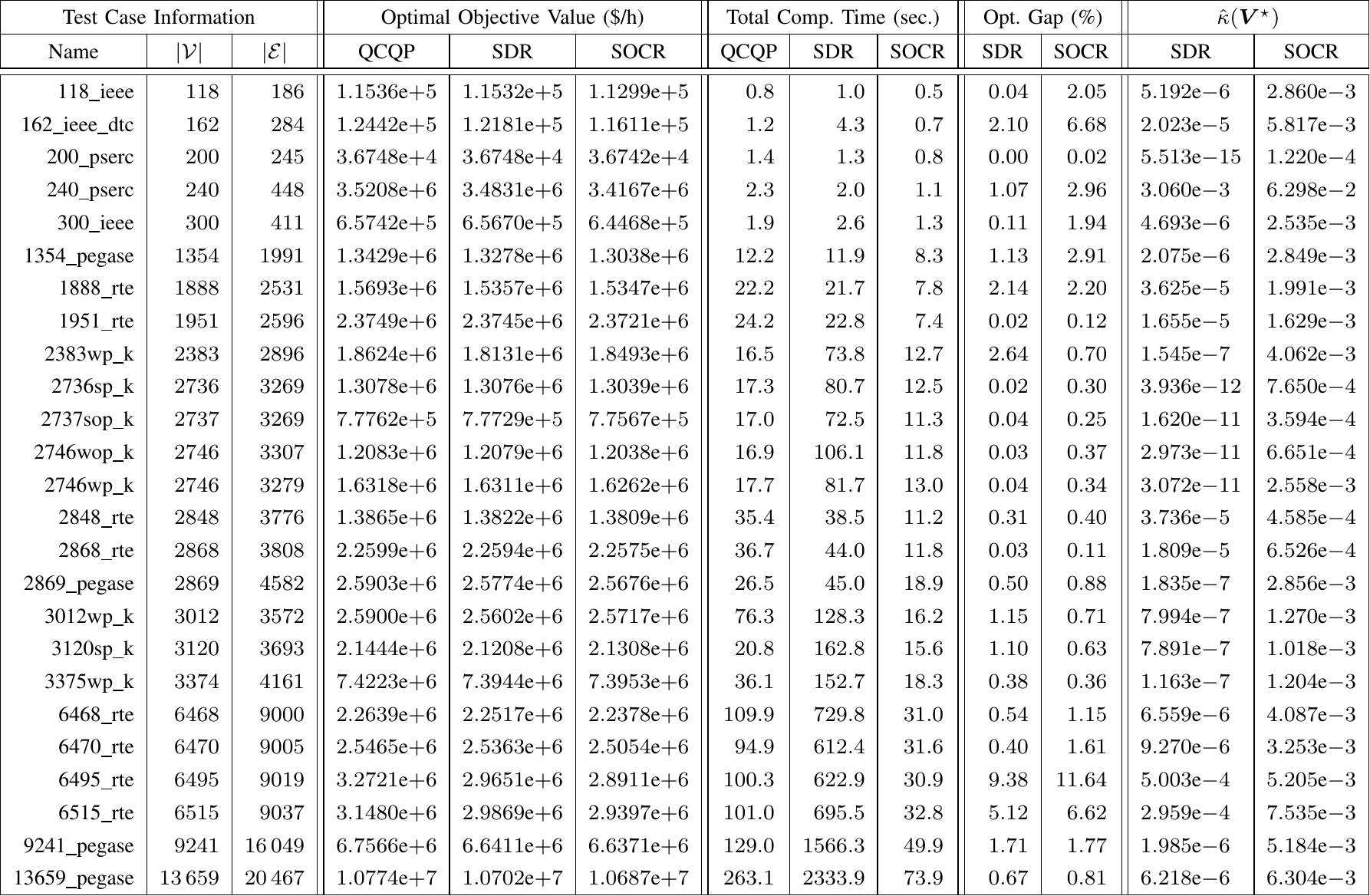}%
\end{table*}

\section{Software Design}
\label{sec:implementation}

This section briefly describes the fundamental design of \hynet as depicted in Fig.~\ref{fig:framework} and relates it to the presented system model and OPF formulations. An OPF study is initiated by loading a scenario from a grid database into a scenario object, where the latter organizes the data using \emph{pandas}~\cite{McKinney2010a} data frames. After optional adjustments, the scenario is used to create an OPF model object, which represents~\eqref{eqn:opf:nonconvex:unified:compact} for the given scenario and whose base class implements the mathematical model in Section~\ref{sec:model} and~\ref{sec:constraints}. This object serves as a \emph{builder}~\cite{Gamma1994a} for the OPF problem, which is represented by an object of a quadratically-constrained quadratic problem (QCQP). This QCQP is solved via a solver interface, where the underlying implementation may solve the nonconvex QCQP~\eqref{eqn:opf:nonconvex:unified:compact}, its SDR~\eqref{eqn:opf:sdp}, or its SOCR~\eqref{eqn:opf:socp}. The solver returns an object containing the result and solution process information, which is routed through a \emph{factory function}~\cite{Gamma1994a} of the OPF model object to obtain an appropriate representation of the result data.

This object-oriented design offers a transparent structure and data flow, while rendering \hynet amenable to extensions. For example, extensions to additional solvers and other relaxations can be implemented via corresponding solver classes. Furthermore, the problem formulation can be customized by subclassing the OPF model or its base class and, by overriding the respective factory function, the result representation can be adjusted accordingly.\footnote{For example, since \hynet v1.2.0 an extension for the \emph{maximum loadability problem}~\cite{Irisarri1997a} is provided, which illustrates such a problem customization.}

\section{Validation and Case Studies}
\label{sec:validation}

In the following, some OPF studies are presented to showcase and validate \hynet. These studies were conducted with \hynet v1.2.1 in Python~3.7.3 using a standard office notebook\footnote{MacBook Pro (2014 series) with an Intel Core i7 CPU and 8\,GB RAM.} and \hynet's solver interfaces for IPOPT~\cite{Wachter2006a} (v3.12.12 with MUMPS~\cite{Amestoy2001a,Amestoy2006a}), MOSEK~\cite{MOSEK2019a} (v9.0.101), and CPLEX~\cite{IBM2019a} (v12.9.0) for the QCQP, (chordal) SDR, and SOCR, respectively. The results also report the \emph{optimality gap}, which is defined as $100\cdot(1-p_{\text{SDR}\vert\text{SOCR}}/p_\text{QCQP})$, where $p_\text{QCQP}$, $p_\text{SDR}$, and $p_\text{SOCR}$ is the optimal objective value of the QCQP, SDR, and SOCR solution, respectively.

\subsection{IEEE PES Power Grid Lib}

The IEEE PES Power Grid Lib (PGLib) OPF~v17.08~\cite{PGLib2017a} is a benchmark for validating AC OPF algorithms and utilized here to validate \hynet. For this study, the PGLib test cases are imported from the \matpower data format into \hynet's grid database format using the command-line interface of \hynet. During the import, polynomial cost functions are converted to piecewise linear functions by sampling the polynomial equidistantly within the generator's active power limits. This study utilizes 10 sample points per polynomial cost function.

Table~\ref{tab:res:pglib:typ} presents the OPF results for the test cases with typical operating conditions and at least 100 buses.\footnote{3375wp\_k contains an isolated bus, which is removed for this study.} It can be observed that these results are consistent with those of \powermodels and \matpower in~\cite{Coffrin2018a}. The minor deviations are explained by the conversion of the cost functions as well as the use of different solvers and thermal flow limits. For the latter, \hynet features the ampacity constraint~\eqref{eqn:crt:ampacity}, while \powermodels and \matpower implement an apparent power flow limit. In contrast to \matpower, which fails on several test cases~\cite{Coffrin2018a}, \hynet successfully solves all problem instances.

In terms of performance, \hynet compares well to \powermodels and \matpower, considering that the latter's computation times in~\cite{Coffrin2018a} are based on high-performance computing servers while the results here were computed on a standard office notebook. Furthermore, the computation times in~\cite{Coffrin2018a} focus on the solver runtime only, while the computation times of \hynet in Table~\ref{tab:res:pglib:typ} additionally include the data verification, modeling, and result processing as well as the computation of a supply-meets-demand (``copper plate'' model) based initial point for QCQPs and a graph traversal based bus voltage recovery for the SDR and SOCR.

In conclusion, these results validate \hynet's functionality and illustrate its robust and competitive performance.

\subsection{Large-Scale Hybrid AC/DC Power System}

To showcase \hynet for a large-scale hybrid AC/DC power system, the model \texttt{case2383wp\_ha.db} of~\cite{Hotz2019a} is considered, which is a MT-HVDC variant of the hybrid system in~\cite{Hotz2018a} and exhibits the \emph{hybrid architecture}, cf. Section~\ref{sec:relaxation:exactness}. It consists of an AC system with 2383 buses and 2392 branches, 16 back-to-back converters, and a multitude of P2P- and MT-HVDC systems that comprise a total of 786 buses and AC/DC converters as well as 488 branches. Its QCQP, SDR, and SOCR OPF is solved by \hynet in 28.0\,sec., 19.2\,sec., and 12.1\,sec. to an optimal objective value of 1867.1\,k\$/h, 1867.1\,k\$/h, and 1866.9\,k\$/h, respectively. For the SDR and SOCR, the optimality gap is 0.00\% and 0.01\% and $\hat{\kappa}(\aOpt{\mV})$ is $5.998\cdot 10^{-15}$ and $1.356\cdot 10^{-14}$, respectively. These results show that \hynet is readily applicable to extensive hybrid AC/DC power systems.

\section{Conclusion}
\label{sec:conclusion}

This paper introduced \hynet, a Python-based open-source OPF framework for hybrid AC/DC grids with point-to-point and radial multi-terminal HVDC systems. The presented system model, which features a unified modeling of AC and DC subgrids as well as a concise converter model, establishes the foundation for \hynet's ease of use and extensibility. It enables a compact mathematical formulation and, additionally, it induces the focused and transparent data management in \hynet. Furthermore, together with the proposed state space relaxation, it simplifies the derivation and study of convex relaxations of the OPF problem for hybrid AC/DC systems and, at the same time, it renders the resulting unified OPF formulation an appropriate data abstraction layer for direct and relaxation-based solution approaches. This mathematical foundation is embedded in a specialized object-oriented software design, which is not only clearly structured but also avoids the need for an explicit interface for extensions due to its extensibility by design. Finally, the presented case studies validated \hynet's functionality, while showing that it also offers a robust and competitive performance. It is hoped that \hynet can serve as a valuable tool for OPF studies with hybrid AC/DC grids and may see future contributions from the community, e.g., solver interfaces for further relaxation techniques and solver implementations.

\appendices

\section{Proof of Theorem~\ref{thm:dcvoltages}}
\label{apx:thm:dcvoltages}

To begin with, consider the following lemma.
\begin{lemma}\label{lem:dcimag}
Consider any $\vf\in\sF$ and $\vs_j\in\sS_j$, for ${j\in\sI}$. Let $\vv\in\nC^\cV$ satisfy~\eqref{eqn:model:crt:powbal} for the given $\vf$ and $\vs_\sI$. Then, $\imag(V_{\eSrc(k)}\conj V_{\eDst(k)})=0$, $\forall k\in\sEDC$.
\end{lemma}
\begin{IEEEproof}
See Appendix~\ref{apx:lem:dcimag}.
\end{IEEEproof}

Thus, for all $k\in\sEDC$, it holds that
\begin{align}\label{eqn:thm:dcvoltages:phase}
	\arg(V_{\eDst(k)}) - \arg(V_{\eSrc(k)}) = r\pi\,,
	\quad\text{for some}\,\
	r\in\nZ\,.
\end{align}
Furthermore, it follows from \eqref{eqn:crt:angle:dcb}, which implements~\eqref{eqn:crt:angle:realpart}, that
\begin{align}\label{eqn:thm:dcvoltages:sign}
	-\pi/2 \leq \arg(V_{\eDst(k)}) - \arg(V_{\eSrc(k)}) \leq \pi/2\,.
\end{align}
Jointly, \eqref{eqn:thm:dcvoltages:phase} and~\eqref{eqn:thm:dcvoltages:sign} imply
\begin{align}
	\arg(V_{\eSrc(k)}) = \arg(V_{\eDst(k)})\,,
	\quad
	\forall k\in\sEDC\,.
\end{align}

Therefore, the phase of all elements in $\vv$ associated with a certain DC subgrid is equal. As the constraints in~\eqref{eqn:opf:nonconvex:full:crt} and~\eqref{eqn:crt:angle:dcb} are quadratic\footnote{\label{foot:quadratic:phaseshift}For $\vx$ and $\aUBar{\vx} = \vx e^{\iu\varphi}$, it follows $\aUBar{\vx}\herm\mA\aUBar{\vx} = e^{-\iu\varphi}\vx\herm\mA\vx e^{\iu\varphi} = \vx\herm\mA\vx$.} in $\vv$ and the individual constraints only involve elements of $\vv$ that are associated with the same subgrid,\footnote{Note that converters decouple the voltages of different subgrids. W.l.o.g., the buses can be numbered such that the constraint matrices are block-diagonal (see also~\cite[App.~B]{Hotz2016a}), where each block relates to a subgrid. Footnote~\ref{foot:quadratic:phaseshift} applies per block, i.e., a blockwise phase shift does not affect the constraints.} it follows that these constraints are invariant with respect to a common phase shift of all elements in $\vv$ for a certain subgrid. Thus, if $\vv$ satisfies~\eqref{eqn:opf:nonconvex:full:crt} and~\eqref{eqn:crt:angle:dcb}, then $\aUBar{\vv}$ in~\eqref{eqn:dcvoltages:reconst} satisfies~\eqref{eqn:opf:nonconvex:full:crt} and~\eqref{eqn:crt:angle:dcb} with equivalent constraint function values, as its construction only involves a common phase shift in the individual DC subgrids.
Furthermore, this implies that $\fLoss(\aUBar{\vv}\aUBar{\vv}\herm,\vf)=\fLoss(\vv\vv\herm,\vf)$, as $\fLoss$ equals the summation of the left-hand side of~\eqref{eqn:model:crt:powbal:p}.

\section{Proof of Lemma~\ref{lem:dcimag}}
\label{apx:lem:dcimag}

For $n\in\sVDC$, it follows from~\eqref{eqn:model:crt:powbal:q}, Def.~\ref{def:dc:conv:injection}, and Def.~\ref{def:dc:injector:injection} that $\vv\herm\mQ_n\vv=0$ and, with Def.~\ref{def:dc:seriescond} and Def.~\ref{def:dc:noshunt}, that\footnote{Alternatively, \eqref{eqn:dcgrid:reactivebalance} can be derived via the reactive power balance equation $\imag(V_n I_n\conj) = 0$ at DC bus $n\in\sVDC$ using the electrical model in Section~\ref{sec:model:electrical}.}
\begin{align}\label{eqn:dcgrid:reactivebalance}
	\sum_{\mathclap{k\in\aSrc{\sBE}(n)}} \aBar{y_k}\imag(V_{\eDst(k)}\conj V_n)
	\ =\,\
	\sum_{\mathclap{k\in\aDst{\sBE}(n)}} \aBar{y_k}\imag(V_{\eSrc(k)}\conj V_n)\,.
\end{align}
Thus, if $\vv$ is not restricted to $\sU$, there may emerge a circulation of ``artificial reactive power'' in DC subgrids.

Now, consider an individual DC subgrid, which is radial by Def.~\ref{def:dcgrid:radial}. Without loss of generality, assume that it comprises no parallel branches (consider their single branch equivalent), consider one of its buses as the reference bus, and let all its branches point toward this reference bus. At all leaf nodes $n$ of this directed tree graph, \eqref{eqn:dcgrid:reactivebalance} reduces to
\begin{align}\label{eqn:dcgrid:reactivebalance:leafnode}
	\aBar{y_k}\imag(V_{\eDst(k)}\conj V_n) = 0
\end{align}
where $k$ is the \emph{only} DC branch connected to DC bus $n = \eSrc(k)$ due to Def.~\ref{def:dcgrid:radial}, Def.~\ref{def:noselfloops}, and the absence of parallel branches. It follows from Def.~\ref{def:lossydcbranches} that $\aBar{y_k}\neq 0$, thus \eqref{eqn:dcgrid:reactivebalance:leafnode} implies
\begin{align}
	\imag(V_{\eDst(k)}\conj V_{\eSrc(k)}) = 0
	\quad\implies\quad
	\imag(V_{\eSrc(k)}\conj V_{\eDst(k)}) = 0
\end{align}
i.e., the inflow of artificial reactive power arising from DC branch $k$ at its destination bus $\eDst(k)$ is zero. Therefore, all DC branches that are connected to leaf nodes do not contribute to the circulation of artificial reactive power and, thus, can be excluded for this analysis. 
Considering the resulting reduced directed tree graph, the above argument can be repeated until the reduced directed tree graph equals the reference node. This proves by induction that $\imag(V_{\eSrc(k)}\conj V_{\eDst(k)})=0$ for all DC branches $k$ of the considered DC subgrid.

Repetition of this inductive argument for all DC subgrids implies that $\imag(V_{\eSrc(k)}\conj V_{\eDst(k)})=0$, $\forall k\in\sEDC$. Note that the proof can be adapted to an arbitrary directionality of DC branches.

\section*{Acknowledgment}

The authors would like to thank Vincent Bode, Michael Mitterer, Christian Wahl, Yangyang He, and Julia Sistermanns for supporting the implementation and testing of \hynet as well as Dominic Hewes, Mohasha Sampath, and Y. S. Foo Eddy for providing feedback on a draft of the manuscript.

\IEEEtriggeratref{64}

\bibliographystyle{IEEEtran}
\bibliography{hotz_utschick_hynet.bib}

% Generated by IEEEtran.bst, version: 1.14 (2015/08/26)
\begin{thebibliography}{10}
\providecommand{\url}[1]{#1}
\csname url@samestyle\endcsname
\providecommand{\newblock}{\relax}
\providecommand{\bibinfo}[2]{#2}
\providecommand{\BIBentrySTDinterwordspacing}{\spaceskip=0pt\relax}
\providecommand{\BIBentryALTinterwordstretchfactor}{4}
\providecommand{\BIBentryALTinterwordspacing}{\spaceskip=\fontdimen2\font plus
\BIBentryALTinterwordstretchfactor\fontdimen3\font minus
  \fontdimen4\font\relax}
\providecommand{\BIBforeignlanguage}[2]{{%
\expandafter\ifx\csname l@#1\endcsname\relax
\typeout{** WARNING: IEEEtran.bst: No hyphenation pattern has been}%
\typeout{** loaded for the language `#1'. Using the pattern for}%
\typeout{** the default language instead.}%
\else
\language=\csname l@#1\endcsname
\fi
#2}}
\providecommand{\BIBdecl}{\relax}
\BIBdecl

\bibitem{European-Commission2017a}
{European Commission}, \emph{The {S}trategic {E}nergy {T}echnology Plan}.\hskip
  1em plus 0.5em minus 0.4em\relax Publications Office of the European Union,
  2017.

\bibitem{Arcia-Garibaldi2018a}
G.~Arcia-Garibaldi, P.~Cruz-Romero, and A.~G\'{o}mez-Exp\'{o}sito, ``Future
  power transmission: {V}isions, technologies and challenges,'' \emph{Renewable
  and Sustainable Energy Reviews}, vol.~94, pp. 285 -- 301, 2018.

\bibitem{IEA2016a}
{OECD/IEA}, \emph{{L}arge-{S}cale {E}lectricity {I}nterconnection: {T}echnology
  and prospects for cross-regional networks}.\hskip 1em plus 0.5em minus
  0.4em\relax Int. Energy Agency, 2016.

\bibitem{Buigues2017a}
G.~Buigues, V.~Valverde, A.~Etxegarai, P.~Egu{\'\i}a, and E.~Torres, ``Present
  and future multiterminal {HVDC} systems: {C}urrent status and forthcoming
  developments,'' in \emph{Proc. Int. Conf. Renewable Energies and Power
  Quality}, Malaga, Spain, Apr. 2017.

\bibitem{Zimmerman2011a}
R.~D. Zimmerman, C.~E. Murillo-S\'{a}nchez, and R.~J. Thomas, ``{MATPOWER}:
  Steady-state operations, planning, and analysis tools for power systems
  research and education,'' \emph{IEEE Trans. Power Syst.}, vol.~26, no.~1, pp.
  12--19, Feb. 2011.

\bibitem{Zimmerman2016a}
\BIBentryALTinterwordspacing
R.~D. Zimmerman, C.~E. Murillo-S\'{a}nchez \emph{et~al.}, ``{MATPOWER}: A
  {Matlab} power system simulation package (v6.0b2),'' Nov. 2016. [Online].
  Available: \url{http://www.pserc.cornell.edu/matpower/}
\BIBentrySTDinterwordspacing

\bibitem{Lincoln2017a}
\BIBentryALTinterwordspacing
R.~Lincoln, ``{PYPOWER}: A port of {MATPOWER} to {P}ython (v5.1.2),'' Jun.
  2017. [Online]. Available: \url{http://github.com/rwl/PYPOWER}
\BIBentrySTDinterwordspacing

\bibitem{Milano2005a}
F.~Milano, ``An open source power system analysis toolbox,'' \emph{IEEE Trans.
  Power Syst.}, vol.~20, no.~3, pp. 1199--1206, Aug. 2005.

\bibitem{Milano2016a}
\BIBentryALTinterwordspacing
------, ``{PSAT}: Matlab-based power system analysis toolbox (v2.1.10),'' Jun.
  2016. [Online]. Available: \url{http://faraday1.ucd.ie/psat.html}
\BIBentrySTDinterwordspacing

\bibitem{Coffrin2018a}
C.~Coffrin, R.~Bent, K.~Sundar, Y.~Ng, and M.~Lubin, ``{PowerModels.jl}: {A}n
  open-source framework for exploring power flow formulations,'' in \emph{Proc.
  Power Systems Computation Conf. (PSCC)}, Jun. 2018.

\bibitem{Coffrin2018b}
\BIBentryALTinterwordspacing
C.~Coffrin \emph{et~al.}, ``{PowerModels.jl}: {A} {J}ulia/{JuMP} package for
  power network optimization (v0.8.5),'' Oct. 2018. [Online]. Available:
  \url{http://github.com/lanl-ansi/PowerModels.jl}
\BIBentrySTDinterwordspacing

\bibitem{Thurner2018a}
L.~Thurner, A.~Scheidler, F.~Sch\"afer, J.~Menke, J.~Dollichon, F.~Meier,
  S.~Meinecke, and M.~Braun, ``Pandapower---{A}n open-source {P}ython tool for
  convenient modeling, analysis, and optimization of electric power systems,''
  \emph{IEEE Trans. Power Syst.}, vol.~33, no.~6, pp. 6510--6521, Nov. 2018.

\bibitem{Thurner2018b}
\BIBentryALTinterwordspacing
------, ``{pandapower}: {A}n easy to use open source tool for power system
  modeling, analysis and optimization with a high degree of automation
  (v1.6.0),'' Sep. 2018. [Online]. Available: \url{http://www.pandapower.org/}
\BIBentrySTDinterwordspacing

\bibitem{Beerten2015a}
J.~Beerten and R.~Belmans, ``Development of an open source power flow software
  for high voltage direct current grids and hybrid {AC}/{DC} systems:
  {MATACDC},'' \emph{IET Generation, Transmission \& Distribution}, vol.~9,
  no.~10, pp. 966--974, 2015.

\bibitem{Wiget2012a}
R.~Wiget and G.~Andersson, ``Optimal power flow for combined {AC} and
  multi-terminal {HVDC} grids based on {VSC} converters,'' in \emph{Proc. IEEE
  PES General Meeting}, Jul. 2012.

\bibitem{Cao2013a}
J.~Cao, W.~Du, H.~F. Wang, and S.~Q. Bu, ``Minimization of transmission loss in
  meshed {AC}/{DC} grids with {VSC}-{MTDC} networks,'' \emph{IEEE Trans. Power
  Syst.}, vol.~28, no.~3, pp. 3047--3055, Aug. 2013.

\bibitem{Feng2014a}
W.~Feng, A.~L. Tuan, L.~B. Tjernberg, A.~Mannikoff, and A.~Bergman, ``A new
  approach for benefit evaluation of multiterminal {VSC}-{HVDC} using a
  proposed mixed {AC}/{DC} optimal power flow,'' \emph{IEEE Trans. Power
  Delivery}, vol.~29, no.~1, pp. 432--443, Feb. 2014.

\bibitem{Zhao2017a}
Q.~Zhao, J.~Garc\'{i}a-Gonz\'{a}lez, O.~Gomis-Bellmunt, E.~Prieto-Araujo, and
  F.~M. Echavarren, ``Impact of converter losses on the optimal power flow
  solution of hybrid networks based on {VSC}-{MTDC},'' \emph{Electric Power
  Systems Research}, vol. 151, pp. 395 -- 403, 2017.

\bibitem{Baradar2013a}
M.~Baradar, M.~R. Hesamzadeh, and M.~Ghandhari, ``Second-order cone programming
  for optimal power flow in {VSC}-type {AC}-{DC} grids,'' \emph{IEEE Trans.
  Power Syst.}, vol.~28, no.~4, pp. 4282--4291, Nov. 2013.

\bibitem{Bahrami2017a}
S.~Bahrami, F.~Therrien, V.~W.~S. Wong, and J.~Jatskevich, ``Semidefinite
  relaxation of optimal power flow for {AC}--{DC} grids,'' \emph{IEEE Trans.
  Power Syst.}, vol.~32, no.~1, pp. 289--304, Jan. 2017.

\bibitem{Beerten2012a}
J.~Beerten, S.~Cole, and R.~Belmans, ``Generalized steady-state {VSC} {MTDC}
  model for sequential {AC}/{DC} power flow algorithms,'' \emph{IEEE Trans.
  Power Syst.}, vol.~27, no.~2, pp. 821--829, May 2012.

\bibitem{Cigre2014b}
{CIGRE Working Group B4.57}, ``Guide for the development of models for {HVDC}
  converters in a {HVDC} grid,'' \emph{CIGRE Brochure 604}, Dec. 2014.

\bibitem{Baradar2013b}
M.~Baradar and M.~Ghandhari, ``A multi-option unified power flow approach for
  hybrid {AC}/{DC} grids incorporating multi-terminal {VSC}-{HVDC},''
  \emph{IEEE Trans. Power Syst.}, vol.~28, no.~3, pp. 2376--2383, Aug. 2013.

\bibitem{Daelemans2008a}
G.~Daelemans, ``{VSC} {HVDC} in meshed networks,'' Master's thesis, Katholieke
  Universiteit Leuven, Leuven, Belgium, 2008.

\bibitem{Wood1996a}
A.~J. Wood and B.~F. Wollenberg, \emph{Power Generation, Operation, and
  Control}, 2nd~ed.\hskip 1em plus 0.5em minus 0.4em\relax Wiley, 1996.

\bibitem{Jabr2006a}
R.~A. Jabr, ``Radial distribution load flow using conic programming,''
  \emph{IEEE Trans. Power Syst.}, vol.~21, no.~3, pp. 1458--1459, Aug. 2006.

\bibitem{Bai2008a}
X.~Bai, H.~Wei, K.~Fujisawa, and Y.~Wang, ``Semidefinite programming for
  optimal power flow problems,'' \emph{International Journal of Electrical
  Power \& Energy Systems}, vol.~30, no. 6-7, pp. 383 -- 392, 2008.

\bibitem{Low2014a}
S.~H. Low, ``Convex relaxation of optimal power flow -- {P}art~{I}:
  Formulations and equivalence,'' \emph{IEEE Trans. Control of Netw. Syst.},
  vol.~1, no.~1, pp. 15--27, Mar. 2014.

\bibitem{Taylor2015a}
J.~A. Taylor, \emph{Convex Optimization of Power Systems}.\hskip 1em plus 0.5em
  minus 0.4em\relax Cambridge, UK: Cambridge University Press, 2015.

\bibitem{Molzahn2015a}
D.~K. Molzahn and I.~A. Hiskens, ``Sparsity-exploiting moment-based relaxations
  of the optimal power flow problem,'' \emph{IEEE Trans. Power Syst.}, vol.~30,
  no.~6, pp. 3168--3180, Nov. 2015.

\bibitem{Coffrin2016a}
C.~Coffrin, H.~L. Hijazi, and P.~V. Hentenryck, ``The {QC} relaxation: A
  theoretical and computational study on optimal power flow,'' \emph{IEEE
  Trans. Power Syst.}, vol.~31, no.~4, pp. 3008--3018, Jul. 2016.

\bibitem{Kocuk2016b}
B.~Kocuk, S.~S. Dey, and X.~A. Sun, ``Strong {SOCP} relaxations for the optimal
  power flow problem,'' \emph{Operations Research}, vol.~64, no.~6, pp.
  1177--1196, Nov. 2016.

\bibitem{Bingane2018a}
C.~Bingane, M.~F. Anjos, and S.~Le~Digabel, ``Tight-and-cheap conic relaxation
  for the {AC} optimal power flow problem,'' \emph{IEEE Trans. Power Syst.},
  vol.~33, no.~6, pp. 7181--7188, Nov. 2018.

\bibitem{Hotz2016a}
M.~Hotz and W.~Utschick, ``A hybrid transmission grid architecture enabling
  efficient optimal power flow,'' \emph{IEEE Trans. Power Syst.}, vol.~31,
  no.~6, pp. 4504--4516, Nov. 2016.

\bibitem{Hotz2018a}
------, ``The hybrid transmission grid architecture: Benefits in nodal
  pricing,'' \emph{IEEE Trans. Power Syst.}, vol.~33, no.~2, pp. 1431--1442,
  Mar. 2018.

\bibitem{Bahrami2018a}
S.~{Bahrami} and V.~W.~S. {Wong}, ``Security-constrained unit commitment for
  {AC}-{DC} grids with generation and load uncertainty,'' \emph{IEEE Trans.
  Power Syst.}, vol.~33, no.~3, pp. 2717--2732, May 2018.

\bibitem{Venzke2019a}
A.~{Venzke} and S.~{Chatzivasileiadis}, ``Convex relaxations of probabilistic
  {AC} optimal power flow for interconnected {AC} and {HVDC} grids,''
  \emph{IEEE Trans. Power Syst.}, vol.~34, no.~4, pp. 2706--2718, Jul. 2019.

\bibitem{Hotz2018b}
\BIBentryALTinterwordspacing
M.~Hotz \emph{et~al.}, ``{\emph{hynet:}} {A}n optimal power flow framework for
  hybrid {AC}/{DC} power systems (v1.0.0),'' Nov. 2018. [Online]. Available:
  \url{http://gitlab.com/tum-msv/hynet}
\BIBentrySTDinterwordspacing

\bibitem{Python-Core-Team2018a}
\BIBentryALTinterwordspacing
{Python Core Team}, ``{Python}: A dynamic, open source programming language
  (v3.7.1),'' Oct. 2018. [Online]. Available: \url{http://www.python.org}
\BIBentrySTDinterwordspacing

\bibitem{Ergun2019a}
H.~{Ergun}, J.~{Dave}, D.~{Van Hertem}, and F.~{Geth}, ``Optimal power flow for
  {AC}--{DC} grids: Formulation, convex relaxation, linear approximation, and
  implementation,'' \emph{IEEE Trans. Power Syst.}, vol.~34, no.~4, pp.
  2980--2990, Jul. 2019.

\bibitem{Bondy1976a}
J.~A. Bondy and U.~S.~R. Murty, \emph{Graph Theory with Applications}.\hskip
  1em plus 0.5em minus 0.4em\relax North Holland, 1976.

\bibitem{Davies2009a}
M.~Davies, M.~Dommaschk, J.~Dorn, J.~Lang, D.~Retzmann, and D.~Soerangr,
  ``{HVDC} {PLUS} -- {B}asics and principle of operation,'' Siemens AG,
  Erlangen, Germany, Tech. Rep., 2009.

\bibitem{Kundur1994a}
P.~Kundur, \emph{Power System Stability and Control}, ser. EPRI Power System
  Engineering Series.\hskip 1em plus 0.5em minus 0.4em\relax McGraw-Hill, 1994.

\bibitem{Bergen2000a}
A.~R. Bergen and V.~Vittal, \emph{Power Systems Analysis}, 2nd~ed.\hskip 1em
  plus 0.5em minus 0.4em\relax Upper Saddle River, NJ, USA: Prentice Hall,
  2000.

\bibitem{Boyd2004a}
S.~Boyd and L.~Vandenberghe, \emph{Convex Optimization}.\hskip 1em plus 0.5em
  minus 0.4em\relax Cambridge, UK: Cambridge University Press, 2004.

\bibitem{Fukuda2001a}
M.~Fukuda, M.~Kojima, K.~Murota, and K.~Nakata, ``Exploiting sparsity in
  semidefinite programming via matrix completion~{I}: General framework,''
  \emph{SIAM J. Optimization}, vol.~11, no.~3, pp. 647--674, 2001.

\bibitem{Jabr2012a}
R.~A. Jabr, ``Exploiting sparsity in {SDP} relaxations of the {OPF} problem,''
  \emph{IEEE Trans. Power Syst.}, vol.~27, no.~2, pp. 1138--1139, May 2012.

\bibitem{Molzahn2013a}
D.~K. Molzahn, J.~T. Holzer, B.~C. Lesieutre, and C.~L. DeMarco,
  ``Implementation of a large-scale optimal power flow solver based on
  semidefinite programming,'' \emph{IEEE Trans. Power Syst.}, vol.~28, no.~4,
  pp. 3987--3998, Nov. 2013.

\bibitem{Andersen2014a}
M.~S. Andersen, A.~Hansson, and L.~Vandenberghe, ``Reduced-complexity
  semidefinite relaxations of optimal power flow problems,'' \emph{IEEE Trans.
  Power Syst.}, vol.~29, no.~4, pp. 1855--1863, Jul. 2014.

\bibitem{Eltved2019a}
A.~Eltved, J.~Dahl, and M.~S. Andersen, ``On the robustness and scalability of
  semidefinite relaxation for optimal power flow problems,'' \emph{Optimization
  and Engineering}, Mar. 2019.

\bibitem{Kim2003a}
S.~Kim and M.~Kojima, ``Exact solutions of some nonconvex quadratic
  optimization problems via {SDP} and {SOCP} relaxations,'' \emph{Computational
  Optimization and Applications}, vol.~26, no.~2, pp. 143--154, Nov. 2003.

\bibitem{Bose2012b}
S.~Bose, S.~H. Low, and K.~M. Chandy, ``Equivalence of branch flow and bus
  injection models,'' in \emph{Proc. 50th Annu. Allerton Conf. Communication,
  Control, and Computing}, Oct. 2012, pp. 1893--1899.

\bibitem{Li2008a}
H.~Li and T.~Adal{\i}, ``Complex-valued adaptive signal processing using
  nonlinear functions,'' \emph{EURASIP J. Advances in Signal Process.}, vol.
  2008, Feb. 2008.

\bibitem{Bazaraa2006a}
M.~S. Bazaraa, H.~D. Sherali, and C.~M. Shetty, \emph{Nonlinear Programming:
  Theory and Algorithms}, 3rd~ed.\hskip 1em plus 0.5em minus 0.4em\relax
  Hoboken, NJ, USA: John Wiley \& Sons, Inc., 2006.

\bibitem{Kirschen2004a}
D.~Kirschen and G.~Strbac, \emph{Fundamentals of Power System Economics}.\hskip
  1em plus 0.5em minus 0.4em\relax John Wiley \& Sons, Ltd, 2004.

\bibitem{Overbye2004a}
T.~J. Overbye, X.~Cheng, and Y.~Sun, ``A comparison of the {AC} and {DC} power
  flow models for {LMP} calculations,'' in \emph{Proc. 37th Annual Hawaii Int.
  Conf. System Sciences}, Jan. 2004, pp. 1--9.

\bibitem{Wang2007a}
H.~Wang, C.~E. Murillo-Sanchez, R.~D. Zimmerman, and R.~J. Thomas, ``On
  computational issues of market-based optimal power flow,'' \emph{IEEE Trans.
  Power Syst.}, vol.~22, no.~3, pp. 1185--1193, Aug. 2007.

\bibitem{Lavaei2012a}
J.~Lavaei and S.~H. Low, ``Zero duality gap in optimal power flow problem,''
  \emph{IEEE Trans. Power Syst.}, vol.~27, no.~1, pp. 92--107, Feb. 2012.

\bibitem{Zhang2013a}
B.~Zhang and D.~Tse, ``Geometry of injection regions of power networks,''
  \emph{IEEE Trans. Power Syst.}, vol.~28, no.~2, pp. 788--797, May 2013.

\bibitem{Farivar2013a}
M.~Farivar and S.~H. Low, ``Branch flow model: Relaxations and convexification
  -- {P}art~{I},'' \emph{IEEE Trans. Power Syst.}, vol.~28, no.~3, pp.
  2554--2564, Aug. 2013.

\bibitem{Low2014b}
S.~H. Low, ``Convex relaxation of optimal power flow -- {P}art~{II}:
  Exactness,'' \emph{IEEE Trans. Control of Netw. Syst.}, vol.~1, no.~2, pp.
  177--189, Jun. 2014.

\bibitem{Madani2015a}
R.~Madani, S.~Sojoudi, and J.~Lavaei, ``Convex relaxation for optimal power
  flow problem: Mesh networks,'' \emph{IEEE Trans. Power Syst.}, vol.~30,
  no.~1, pp. 199--211, Jan. 2015.

\bibitem{Nick2018a}
M.~Nick, R.~Cherkaoui, J.~L. Boudec, and M.~Paolone, ``An exact convex
  formulation of the optimal power flow in radial distribution networks
  including transverse components,'' \emph{IEEE Trans. Automatic Control},
  vol.~63, no.~3, pp. 682--697, Mar. 2018.

\bibitem{Gan2014b}
L.~Gan and S.~H. Low, ``Optimal power flow in direct current networks,''
  \emph{IEEE Trans. Power Syst.}, vol.~29, no.~6, pp. 2892--2904, Nov. 2014.

\bibitem{Tan2015a}
C.~W. Tan, D.~W.~H. Cai, and X.~Lou, ``Resistive network optimal power flow:
  Uniqueness and algorithms,'' \emph{IEEE Trans. Power Syst.}, vol.~30, no.~1,
  pp. 263--273, Jan. 2015.

\bibitem{McKinney2010a}
W.~McKinney, ``Data structures for statistical computing in {P}ython,'' in
  \emph{Proc. 9th Python in Science Conf.}, 2010, pp. 51 -- 56.

\bibitem{Gamma1994a}
E.~Gamma, R.~Helm, R.~Johnson, and J.~Vlissides, \emph{Design Patterns:
  Elements of Reusable Object-Oriented Software}, ser. Professional Computing
  Series.\hskip 1em plus 0.5em minus 0.4em\relax Addison-Wesley, 1994.

\bibitem{Irisarri1997a}
G.~D. {Irisarri}, X.~{Wang}, J.~{Tong}, and S.~{Mokhtari}, ``Maximum
  loadability of power systems using interior point nonlinear optimization
  method,'' \emph{IEEE Trans. Power Syst.}, vol.~12, no.~1, pp. 162--172, Feb.
  1997.

\bibitem{Wachter2006a}
A.~W{\"a}chter and L.~T. Biegler, ``On the implementation of an interior-point
  filter line-search algorithm for large-scale nonlinear programming,''
  \emph{Mathematical Programming}, vol. 106, no.~1, pp. 25--57, Mar. 2006.

\bibitem{Amestoy2001a}
P.~R. Amestoy, I.~S. Duff, J.~Koster, and J.-Y. L'Excellent, ``A fully
  asynchronous multifrontal solver using distributed dynamic scheduling,''
  \emph{SIAM J. Matrix Analysis and Applications}, vol.~23, no.~1, pp. 15--41,
  2001.

\bibitem{Amestoy2006a}
P.~R. Amestoy, A.~Guermouche, J.-Y. L'Excellent, and S.~Pralet, ``Hybrid
  scheduling for the parallel solution of linear systems,'' \emph{Parallel
  Computing}, vol.~32, no.~2, pp. 136--156, 2006.

\bibitem{MOSEK2019a}
\BIBentryALTinterwordspacing
{MOSEK ApS}, \emph{The {MOSEK} Optimizer {API} for {P}ython manual.
  {V}ersion~9.0.}, 2019. [Online]. Available:
  \url{http://docs.mosek.com/9.0/pythonapi.pdf}
\BIBentrySTDinterwordspacing

\bibitem{IBM2019a}
\BIBentryALTinterwordspacing
{IBM}, \emph{{IBM} {ILOG} {CPLEX} Optimization Studio (v12.9)}, 2019. [Online].
  Available: \url{http://www.ibm.com/products/ilog-cplex-optimization-studio}
\BIBentrySTDinterwordspacing

\bibitem{PGLib2017a}
\BIBentryALTinterwordspacing
{IEEE PES Task Force on Benchmarks for Validation of Emerging Power System
  Algorithms}, \emph{{P}ower {G}rid {L}ib: {B}enchmarks for Validating Power
  System Algorithms (v17.08)}, Aug. 2017. [Online]. Available:
  \url{http://github.com/power-grid-lib/pglib-opf}
\BIBentrySTDinterwordspacing

\bibitem{Hotz2019a}
\BIBentryALTinterwordspacing
M.~{Hotz (Curator)}, ``{\emph{hynet}} {G}rid {D}atabase {L}ibrary (v1.4),''
  Aug. 2019. [Online]. Available:
  \url{http://gitlab.com/tum-msv/hynet-databases}
\BIBentrySTDinterwordspacing

\end{thebibliography}

\vfill

\end{document}